\documentclass[12pt,reqno]{amsart}
\usepackage{amsmath}
\usepackage{amssymb}
\usepackage{amsthm}
\usepackage{mathrsfs}

\usepackage{latexsym}
\usepackage{amssymb}

\newtheorem{theorem}{Theorem}[section]
\newtheorem{lemma}{Lemma}[section]
\newtheorem{corollary}{Corollary}[section]
\newtheorem{remark}{Remark}[section]
\newtheorem{definition}{Definition}[section]
\newtheorem{proposition}{Proposition}[section]
\newtheorem{example}{Example}[section]
\newtheorem{assumption}{Assumption}[section]
\numberwithin{equation}{section}
\newcommand{\bth}{\begin{theorem}}
\newcommand{\ethe}{\end{theorem}}

\newcommand{\bre}{\begin{remark}}
\newcommand{\ere}{\end{remark}}

\newcommand{\ble}{\begin{lemma}}
\newcommand{\ele}{\end{lemma}}
\newcommand{\bde}{\begin{definition}}
\newcommand{\ede}{\end{definition}}

\newcommand{\bco}{\begin{corollary}}
\newcommand{\eco}{\end{corollary}}

\newcommand{\bpr}{\begin{proposition}}
\newcommand{\epr}{\end{proposition}}

\newcommand{\bexer}{\begin{exercise}}
\newcommand{\eexer}{\end{exercise}}
\newcommand{\breh}{\begin{hint}}
\newcommand{\ereh}{\end{hint}}

\newcommand{\halmos}{\hfill \qed}

\newcommand{\bexam}{\begin{example}}
\newcommand{\eexam}{\end{example}}

\newcommand{\pr} {{\bf Proof.}}

\newcommand{\bfi}{\begin{fig}}
\newcommand{\efi}{\end{fig}}

\newcommand{\beao}{\begin{eqnarray*}}
\newcommand{\eeao}{\end{eqnarray*}\noindent}

\newcommand{\beam}{\begin{eqnarray}}
\newcommand{\eeam}{\end{eqnarray}\noindent}
\newcommand{\E}{\mathbf{E}}
\newcommand{\PP}{\mathbf{P}}

\newcommand{\xto}{x\to\infty}

\newcommand{\bA}{\overline{A}}
\newcommand{\bH}{\overline{H}}
\newcommand{\bF}{\overline{F}}

\newcommand{\bG}{\overline{G}}
\newcommand{\bV}{\overline{V}}

\newcommand{\bbr}{{\mathbb R}}

\newcommand{\bbb}{{\mathbb B}}

\newcommand{\bbn}{{\mathbb N}}

\newcommand{\vep}{\varepsilon}

\begin{document}

\title[Subexponentiality and interplay of insurance and financial risks]{Multivariate subexponentiality and interplay of insurance and financial risks in a renewal risk model}
\author[D.G. Konstantinides, C. D. Passalidis]{Dimitrios G. Konstantinides, Charalampos  D. Passalidis}

\address{Dept. of Statistics and Actuarial-Financial Mathematics,
University of the Aegean,
Karlovassi, GR-83 200 Samos, Greece}
\email{konstant@aegean.gr,\;sasd24009@sas.aegean.gr.}
\begin{abstract}
During last decade, although the study of continuous time risk models with heavy-tailed claims and presence of stochastic returns became quite popular, due to its practical applications, still it is lacking the emphasis in modeling the dependence among these two fundamental risks, that seems necessary for the modern insurance industry. In this paper we consider a multivariate risk model with common renewal process, while the logarithmic returns of the insurer's investment portfolio, are described by a L\'{e}vy process.

In the two main results are established an asymptotic expression for the entrance probability of the discounted aggregate claims in some 'rare' sets $x\,A$. This asymptotic expression highlights the multivariate linear single big jump principle in asymptotic behavior of these probabilities.

In the first result, we are restricted in the case where the insurer makes risk-free investments, and hence we consider a non-negative L\'{e}vy process. We assume that the claim vectors follow a distribution from the class $\mathcal{A}_A^*$, introduced here, and represents a negligibly smaller subclass of class $\mathcal{S}_A$ of multivariate subexponential distributions on $A$, since the additional requirement for positive lower Karamata index, looks as a quit mild condition. Further, we consider that the insurance and financial risks, satisfy a weak, but very general dependence structure.

In the second result, we allow arbitrarily dependence between the two risks, and we 
assume that the distribution of their product, at each renewal epoch, belongs to the 
class $(\mathcal{D} \cap \mathcal{A})_A \subsetneq \mathcal{S}_A$. In this theorem we also permit risky-investment, putting a condition to Laplace exponent of the L\'{e}vy process, related with the upper Matuszewska index of the distribution for the product insurance and financial risks. We also note that even in the special one-dimensional sub-case the main results are new.

Furthermore, we present two examples, where we demand only conditions for the marginal distributions of both risks and their dependence structure. In the first example, we consider the weak dependence structure, used in the first theorem, and the insurance risks dominate against the financial risks. In the second example we take a strong dependence structure, and we also allow the cases where the financial risks dominate against the insurance risks, or the two risks have equal "heaviness". Both examples, under the restriction on multivariate regularly varying distributions provide more explicit and elegant relations in relation with that established in the main results. Finally, we give an application of the results for the asymptotic estimation of the 
infinite-time ruin probability.
\end{abstract}

\maketitle
\vspace{3mm}
\vspace{3mm}



\section{Introduction} \label{sec.KLP.1}

The modern insurance companies handle more than one risky portfolios, and in order to be viable and the same time competitive, invest their surplus in risk-free and in risky assets. Thus, the insurer faces two types of risks, the insurance risks, caused by insurance claims, and the financial risks which caused by risky investments. We find an increasing number of papers in literature on this topic, that study multidimensional risk models with heavy-tailed insurance and financial risks, see for example \cite{li:2016}, \cite{yang:su:2023}, \cite{konstantinides:passalidis:2024j}, \cite{chen:konstantinides:passalidis:2025}. 
However, in these papers on one hand the insurance risks follow some heavy-tailed distribution and dominate on the financial risks, namely the insurance risks have heavier tail than the financial risks, and on the other hand the insurance and financial risks are independent. With respect to these two drawbacks, we find similar difficulties in one dimensional risk models too, with only few exceptions as for example in \cite{guo:2022}, \cite{yang:fan:yuen:2023} in one dimensional risk models and finite-time horizon results, and \cite{cheng:konstantinides:wang:2024} for multidimensional risk model and infinite-time horizon results. 

So, in this paper we consider an insurer, who operates $d$-lines of business, with $d \in \bbn$, and the claims are described by the sequence of vectors $\{{\bf X}^{(i)}\,,\;i\in \bbn\}$, representing independent and identically distributed (symbolically, i.i.d.) copies of the non-negative random vector ${\bf X}$ with distribution $F$. Let us note, that each claim vector ${\bf X}^{(i)}=(X_{1}^{(i)},\,\ldots,\,X_{d}^{(i)})$, with $i \in \bbn$, can have zero components, but not all of them, and arrives at the time moment $\tau_i$, $i \in \bbn$. Thus, we suppose that the sequence of $\{\tau_i\,,\;i\in N\}$, with $\tau_0=0$, represents a (renewal) counting process $N(t) :=\sup\{i\in \bbn\;:\;\tau_i \leq t\}\,,\,\, t\geq 0$,  for any $t\geq 0$, with $\sup \emptyset = 0$ conventionally, and with finite renewal function
\beao
\lambda(t) := \E[N(t)] = \sum_{i=1}^{\infty} \PP(\tau_i \leq t)\,.
\eeao
In order to avoid the trivial cases, in all over the paper the asymptotic relations holds on $\Lambda$, with $\Lambda:=\{t\;:\;\lambda(t) > 0 \}$. Further, we assume that the insurer invests his surplus in risk-free and in risky assets and the logarithmic return process of the investment portfolio is described by a L\'{e}vy process $\{R(t)\,,\;t\geq 0\}$, which starts from zero, with Laplace exponent of the form $\phi(s) = \log\E\left[e^{-s\,R(1)}\right]$, for any $s \in \bbr$, where if $\phi(s)< \infty$, then we obtain that
\beam \label{eq.KLP.1.1}
\E\left[e^{-s\,R(t)}\right] = e^{t\,\phi(s)}\,, 
\eeam
for any $t\geq 0$, see in \cite[Prop. 3.14]{cont:tankov:2004}. This way, for some fixed $T \in\Lambda$, the discounted aggregate claims of insurer are presented as 
\beam \label{eq.KLP.1.2}
{\bf D}(T) =\sum_{i=1}^{N(T)} {\bf X}^{(i)}\,e^{-R(\tau_i)}=\left( 
\begin{array}{c}
\sum_{i=1}^{N(T)} X_{1}^{(i)}\,e^{-R(\tau_{i}) } \\ 
\vdots \\ 
\sum_{i=1}^{N(T)} X_{d}^{(i)}\,e^{-R(\tau_{i}) } 
\end{array} 
\right)\,,
\eeam 
and in case of infinite time horizon, we have that
\beam \label{eq.KLP.1.3}
{\bf D}(\infty) =\sum_{i=1}^{\infty} {\bf X}^{(i)}\,e^{-R(\tau_i)}=\left( 
\begin{array}{c}
\sum_{i=1}^{\infty} X_{1}^{(i)}\,e^{-R(\tau_{i}) } \\ 
\vdots \\ 
\sum_{i=1}^{\infty} X_{d}^{(i)}\,e^{-R(\tau_{i}) } 
\end{array} 
\right)\,.
\eeam 
We focus our attention on entrance probability of the ${\bf D}(\infty)$ into some 'rare sets' $x\,A$, with $A$ some set from a wide set family $\mathscr{R}$, described in Section 2, as $\xto$. Namely, we are interested to asymptotic expressions of the form
\beam \label{eq.KLP.1.4}
\PP\left[{\bf D}(\infty) \in x\,A\right] \sim \int_0^{\infty} \PP\left[{\bf X}\,e^{-\,R(s)}\in x\,A\right] \,\lambda(ds)\,, 
\eeam
which indicates the presence of the multivariate linear-single big jump principle, for the discounted aggregate claims.

Under various assumptions about the distribution $F$ of the claim vectors and about the investment process, relation \eqref{eq.KLP.1.4}, and its finite-time analogue was established, see for example in \cite{konstantinides:passalidis:2024j}, \cite{chen:konstantinides:passalidis:2025}, \cite{passalidis:2025}. We already mentioned, that in all these papers there exists assumption of domination of the insurance risks versus the financial risks, and independence between these two kinds of risk. As we shall see later, these probabilities are directly related with several kinds of ruin probability.

The rest of this paper is organized as follows. In Section 2, after introduction 
of the necessary preliminaries for the multivariate heavy-tailed distributions, we introduce a new multivariate distribution class. This class is denoted by $\mathcal{A}_A^*$, and it contains all the multivariate subexponential distributions, $\mathcal{S}_A$, with positive lower Karamata index, see Definition \ref{def.KPT.c.1}. For this class, we show that from the one hand side it contains the famous multivariate regular variation, and from the other hand side is slightly smaller than $\mathcal{S}_A$. We also suggest a new dependence structure that depicts the dependence between the insurance and financial risks. This dependence is a general enough weak dependence structure.

In Section 3, we study the asymptotic behavior of the probability of entrance of discounted aggregate claims over infinite time horizon, and we establish relation \eqref{eq.KLP.1.4}, for the case when the claim vectors follow distribution $F$ from the class $\mathcal{A}_A^*$, the investment portfolio logarithmic return process, is described by a non-negative L\'{e}vy process and the dependence between the two risks is given by the weak dependence structure from Section 2. In case we restrict $F$ on the class of multivariate regularly varying distributions, a more direct expression than the formula \eqref{eq.KLP.1.4} is provided. 

In Section 4, we consider arbitrarily dependent insurance and financial  risks and we assume that their product at each renewal epoch has common distribution $H \in (\mathcal{D}\cap\mathcal{A})_A$, see in Subsection 2.1, and we replace the assumption of non-negativity of the L\'{e}vy process with a more relaxed condition on its Laplace exponent. Under these conditions we establish relation again \eqref{eq.KLP.1.4}. Further we give two examples to show the crucial role of dependence and the way how to 
check the assumption, $H\in (\mathcal{D} \cap \mathcal{A})_A$, based only on conditions of the marginal distributions of two risks and their dependence structure. In the first example we employ the weak dependence of Assumption \ref{ass.KPT.2.2}, and in the other one a strong dependence.

The interesting thing in the second example is that permits also the domination of the financial risks against the insurance risks and also the case, in which these two are balanced, which is not previously examined in continuous time risk models. Finally, in Section 5, we present an applications on the ruin probability over infinite-time horizons.

\section{Preliminaries} \label{sec.KLP.2}

In what follows the limit relations hold as $\xto$, except otherwise stated. For two positive functions $f$ and $g$, we denote by $f(x) \sim c\,g(x)$, for some constant $c>0$, if it holds $\lim f(x)/g(x) = c$, by $f(x)=o[g(x)]$, if it holds $\lim f(x)/g(x)= 0$, by $f(x)=O[g(x)]$, if it holds $\limsup f(x)/g(x) < \infty$. Further, we write $f(x) \asymp g(x)$, if $f(x)=O[g(x)]$ and $g(x)=O[f(x)]$ hold simultaneously. By $f^{\leftarrow}$ we denote the c\'{a}gl\'{a}d inverse of $f$.

In what follows, the $d$-dimensional vectors are denoted by bold script. For any two vectors ${\bf x},\,{\bf y}$ and some positive quantity $\lambda>0$, we define the operations component-wisely as usual, namely ${\bf x}\pm{\bf y}=(x_1\pm y_1,\,\ldots,\,x_d\pm y_d)$ and $\lambda\,{\bf x}=(\lambda\,x_1,\,\ldots,\,\lambda\,x_d)$, respectively, while with ${\bf x}^T$ the inverse of ${\bf x}$. For any set $\bbb \in \bbr^d$ we denote by $\overline{\bbb}$ its closed hull, by $\partial \bbb$ its border, by $\bbb^c$ its complement set, and by ${\bf 1}_{\bbb}$ its indicator function. We say that the set $\bbb$ is increasing, if for any ${\bf x} \in \bbb$ and any ${\bf y} \in  \bbr^d_{+}:=[0,\,\infty)^d$, it holds that ${\bf x}+{\bf y} \in \bbb$. All the random vectors have distribution with support on $\bbr^d_{+}$

For any set $\bbb \in \bbr^d_{+} \setminus \{{\bf 0}\}$, the notation of $d$-dimensional, with $d \in \bbn$, positive functions ${\bf f}$ and ${\bf g}$ remain in tact with one dimensional one, for example we denote ${\bf f}(x\,\bbb)=o[{\bf g}(x\,\bbb)]$, if it holds
\beao
\lim \dfrac {{\bf f}(x\,\bbb)}{{\bf g}(x\,\bbb)} =0\,.
\eeao
For two real numbers $a$ and $b$ we denote $a\vee b := \max\{a,\,b\}$, $a\wedge b := \min\{a,\,b\}$. Finally, we denote by $\Theta \sim V$ when a random variable (or vector) $\Theta$ has distribution $V$, while we denote by $\bV(x)=1-V(x)$ the distribution tail of $V$, for any $x \in \bbr$. Also by $V^{n*}$ is denoted the $n$-fold convolution of $V$ with itself. We denote by $S(V)$ the support of distribution $V$, or we write $S(\Theta)$, when we have in mind the random variable $\Theta\sim V$.

\subsection{Heavy-tailed distributions}

The multivariate subexponential distribution, defined in \cite{samorodnitsky:sun:2016}, is based on the following family of sets:
\beam \label{eq.KLP.2.1}
\mathscr{R}:=\{A \subsetneq \bbr^d\;:\;A\;{\text open},\; {\text increasing},\; A^c \;{\text convex},\; {\bf 0} \notin \bA \}\,.
\eeam
For some random vector ${\bf Z}$ with distribution $V$, by \cite[Lem. 4.5]{samorodnitsky:sun:2016} we find that
\beam \label{eq.KPT.2.2}
Z_A:=\sup\{u\;:\; {\bf Z} \in u\,A \}\,,
\eeam
represents a random variable, with proper distribution $V_A$, with tail of the form
\beam \label{eq.KLP.2.3}
\bV_A(x)=\PP({\bf Z} \in x\,A) = \PP\left[\sup_{{\bf p} \in I_A} {\bf p}^T\,{\bf Z} >x \right]\,,
\eeam
for any $x>0$, where for $I_A \subsetneq \bbr^d$, see in \cite[Lem. 4.3(c)]{samorodnitsky:sun:2016}. So, for some $A \in \mathscr{R}$, the class of multivariate subexponential distributions on $A$, symbolically $\mathcal{S}_A$, was defined through relations \eqref{eq.KLP.2.1} - \eqref{eq.KLP.2.3}. Namely, we say that $V \in \mathcal{S}_A$, if $V_A \in \mathcal{S}$, that means for any (or, equivalently, for some) integer $n\geq 2$, it holds
\beao
\lim \dfrac {\overline{V_A^{n*}}(x)}{\bV_A(x)}= n\,.
\eeao
 Similarly in \cite{konstantinides:passalidis:2024g} we find the following distribution classes. We say that a distribution $V$ belongs to the class of multivariate long tailed distributions on $A$, symbolically $V \in\mathcal{L}_A$, if $V_A \in \mathcal{L}$, namely we have that 
\beao
\lim \dfrac {\overline{V}_A(x-a)}{\bV_A(x)}=1\,,
\eeao
for any (or, equivalently, for some) $a>0$. It is well known that if $V_A \in \mathcal{L}$, then there exists a function $a\;:\;[0,\,\infty) \to (0,\,\infty)$, such that $a(x) \to \infty$, $a(x) = o(x)$ and $\bV_A(x\pm a(x)) \sim \bV_A(x)$. This function called insensitivity function of $V_A$, see \cite[Sect. 2.8]{foss:korshunov:zachary:2013} for more details.

We say that a distribution $V$ belongs to the class of multivariate dominatedly varying distributions on $A$, symbolically $V \in\mathcal{D}_A$, if $V_A \in \mathcal{D}$, namely we have that 
\beao
\limsup \dfrac {\overline{V}_A(b\,x)}{\bV_A(x)}< \infty\,,
\eeao
for any (or, equivalently, for some) $b \in (0,\,1)$. Further we write $V \in (\mathcal{D}\cap \mathcal{L})_A$, if $V_A \in \mathcal{D}\cap \mathcal{L}$, and since by \cite{goldie:1978} we have that $\mathcal{D}\cap \mathcal{L} \equiv \mathcal{D}\cap \mathcal{S}$, it follows that $(\mathcal{D}\cap \mathcal{L})_A\equiv (\mathcal{D}\cap \mathcal{S})_A$.

From \cite{konstantinides:passalidis:2024h} we find the class of multivariate positively decreasing distributions on $A$, symbolically $\mathcal{P_D}_A$ defined as follows. We say  $V \in \mathcal{P_D}_A$, if $V_A \in \mathcal{P_D}$, namely for any (or, equivalently, for some) $v>1$ it holds that
\beao
\limsup \dfrac {\overline{V}_A(v\,x)}{\bV_A(x)}< 1\,.
\eeao

In one dimension case, class $\mathcal{P_D}$ is a large enough distribution family and contains distributions with either heavy or light tails. Concretely, it contains all the distributions with infinite right endpoint in their support, which are NOT  extended slowly varying, see \cite[Sect. 7.4]{konstantinides:2018}. Hence, the class $\mathcal{A}=\mathcal{S}\cap \mathcal{P_D}$, introduced in \cite{konstantinides:tang:tsitsiashvili:2002}, is slightly smaller than $\mathcal{S}$, with many applications in actuarial science, as mentioned by \cite{tang:2006}. From \cite{bardoutsos:konstantinides:2011} we have the study of the classes $\mathcal{D}\cap \mathcal{P_D}$ and $\mathcal{D}\cap \mathcal{A}$. We refer to the reader also in \cite{konstantinides:passalidis:2025b} for an overview of properties of $\mathcal{P_D}$ and related classes. So for $A \in \mathscr{R}$, we have that $V \in \mathcal{B}_A$ if $V_A \in \mathcal{B}$, with $\mathcal{B} \in \{\mathcal{A},\,\mathcal{D}\cap \mathcal{P_D},\,\mathcal{D}\cap \mathcal{A}\}$.

For these previous classes we write $\mathcal{B}_\mathscr{R}:=\bigcap_{A \in \mathscr{R}} \mathcal{B}_A$, with $\mathcal{B} \in \{\mathcal{S},\,\mathcal{L},\,\mathcal{D},\,\mathcal{D}\cap \mathcal{L},\,\mathcal{P_D},\,\mathcal{A},\,\mathcal{D}\cap \mathcal{A}\}$.  

Now, we should remind some indexes that have direct relation with the characterization of distribution classes. For an one dimensional distribution $G$ the upper and lower Matuszewska indexes, are denoted as follows
\beao
J_G^+ :=\inf\left\{- \dfrac{\ln \bG_{*}(v)}{\ln v}\;:\;v>1 \right\}= -\lim_{v \rightarrow \infty }\dfrac{\ln \bG_{*}(v )}{\ln v }\,, \\[2mm]
J_G^- :=\sup \left\{- \dfrac{\ln \bG^{*}(v)}{\ln v}\;:\;v>1 \right\}= -\lim_{v \rightarrow \infty }\dfrac{\ln \bG^{*}(v )}{\ln v }\,.
\eeao
where
\beao
\bG_{*}(v):=\liminf_{x\rightarrow \infty }\dfrac{\bG(v\,x)}{\bG(x)}\,, \qquad \bG^{*}(v):=\limsup_{x\rightarrow \infty }\dfrac{\bG(v\,x)}{\bG(x)}\,,
\eeao
for any $v>0$. It is interesting to notice that for any distribution $G$, with unbounded from above support, we obtain that $0\leq J_G^-\leq J_G^+ \leq \infty$. Furthermore it is well known that $G \in \mathcal{D}$ if and only if $J_G^+ < \infty$, and $G \in \mathcal{P_D}$ if and only if $J_G^- > 0$. In case of regularly varying distribution $G$ with index $\alpha\in (0,\infty)$, symbolically $G \in \mathcal{R}_{-\alpha}$, (remind that  $G \in \mathcal{R}_{-\alpha}$, if it holds $\lim \bar{G}(tx)/\bar{G}(x) =t^{-\alpha}$, for any $t > 0$) then $ J_G^-= J_G^+ = \alpha$. Also, it is well known that if $G \in\mathcal{D}$, then for any $p> J_G^+$, it holds that
\beam \label{eq.KP.2.3a}
x^{-p}=o\left[\bG(x)\right]\,.
\eeam 
For more details about these indexes see \cite[Sect. 2.1.2]{bingham:goldie:teugels:1987} and \cite[Sect. 2.4]{leipus:siaulys:konstantinides:2023}.

Now we present the class  of multivariate regularly varying distributions, symbolically $MRV$, in standard form. Let ${\bf Z}$ a random vector with distribution $V$, then we say that it belongs to $MRV$, if there exists some one-dimensional distribution $G \in \mathcal{R}_{-\alpha}$, for some $\alpha \in (0,\,\infty)$, and some non-degenerate to zero, Radon measure $\mu$, such that it holds
\beao
\lim \dfrac{\PP({\bf Z} \in x\,\bbb)}{\bG(x)} = \mu(\bbb)\,,
\eeao
for any $\mu$-continuous Borel set $\bbb \subsetneq [0,\,\infty]^d \setminus \{{\bf 0}\}$. In this case we write $F \in MRV(\alpha,\,\mu)$. We want to notice that, the last relation is  equivalent with:
\beao
\lim x\,\PP\left[ \dfrac{\bf Z}{U_G(x)} \in \bbb \right] = \mu(\bbb)\,,
\eeao
for any Borel set $\bbb \subseteq [0,\,\infty]^d \setminus \{{\bf 0}\}$, which is $\mu$-continuous. The normalization function $U_G$ is given by the relation
\beao
U_G(x)=\left(\dfrac 1{\bG} \right)^{\leftarrow}(x)\,,
\eeao
Hence, we observe that the measure $\mu$ is homogeneous, namely for any $\bbb \subsetneq [0,\,\infty]^d \setminus \{{\bf 0}\}$ and any real $\lambda>0$ we have that
\beao
\mu(\lambda\,\bbb) = \lambda^{-\alpha}\,\mu(\bbb)\,,
\eeao
while if the condition $\mu(({\bf 1},\,\vec{\infty}])>0$ holds, then this implies the fact that the components of ${\bf Z}=(Z_1,\,\ldots,\,Z_d)$, are asymptotically dependent, and have distributions such that 
\beao
\PP(Z_i >x ) \sim \mu(({\bf 1}_i,\,\infty]) \bG(x)\,,
\eeao 
where $0 <\mu(({\bf 1},\,\vec{\infty}]) \leq \mu(({\bf 1}_i,\,\infty]) < \infty$ 
with ${\bf 1}_i$ representing the vector, whose $i$-th components is unit and all the others are zero. Thus, from the last relation we obtain that the $Z_1,\,\ldots,\,Z_d$ follow again distributions belonging to $ \mathcal{R}_{-\alpha}$. 

In general, $MRV$ has found many applications in actuarial and financial mathematics, see for example in \cite{li:2016}, \cite{li:2022a}, \cite{yang:su:2023}, \cite{chen:cheng:zheng:2025}. Let notice that in most of these papers appears the condition $\mu(({\bf 1},\,\vec{\infty}])>0$. Finally, from \cite[Prop. 3.1]{konstantinides:passalidis:2024h} in combination with the corresponding one-dimensional inclusions, see for example in \cite[Ch. 2]{leipus:siaulys:konstantinides:2023} we obtain that 
\beam \label{eq.KLP.2.10}
MRV\subsetneq (\mathcal{D} \cap \mathcal{A})_\mathscr{R} \subsetneq  (\mathcal{D} \cap \mathcal{L})_\mathscr{R}\subsetneq  \mathcal{S}_\mathscr{R}\subsetneq  \mathcal{L}_\mathscr{R}\,,
\eeam 
where the inclusions in \eqref{eq.KLP.2.10} remain in tact if instead  of $\mathscr{R}$ we put $\mathcal{B}_A$ for some $A \in \mathscr{R}$, with $\mathcal{B} \in \{\mathcal{D} \cap \mathcal{A},\,\mathcal{D} \cap \mathcal{L},\,\mathcal{S},\,\mathcal{L}\}$. In  \eqref{eq.KLP.2.10} as $MRV$ we keep in mind the union of all the $MRV(\alpha,\,\mu)$ for any $\alpha \in (0,\,\infty)$ and any non-degenerate to zero Radon measure $\mu$.

There are many monographs devoted to application of $MRV$, in several other branches of 
applied probability, see for example \cite{resnick:2007}, \cite{buraczewski:damek:mikosch:2016}, \cite{samorodnitsky:2016}, \cite{kulik:soulier:2020}. Unfortunately, in most of them is considered the standard $MRV$, which facilitates the multivariate extreme value theory but does not help to practical applications of actuarial and financial science, for the following reasons:
\begin{enumerate} 
	\item[(i)]
	$MRV$ does NOT cover 'moderately  heavy' tails in its components, as for example in case of Lognormal and Weibull distributions.
	\item[(ii)]
	The standard $MRV$ restricts the inhomogeneity with respect to heaviness of the tails 
	in the portfolios, something that appears frequently, especially in insurance companies with big number of portfolios.
\end{enumerate} 

To face the (i) and (ii), we work in the classes $\mathcal{A}_A^*$ and $(\mathcal{D}\cap \mathcal{A})_A$, for $A \in \mathscr{R}$, since both classes permits cases of inhomogeneous components with respect to heaviness of the tails in the components of ${\bf X}$, see for example in \cite[Exam. 4.1 - 4.4]{konstantinides:liu:passalidis:2025}.

\subsection{On class $\mathcal{A}_A^*$}

In this subsection we introduce the multivariate version of distribution class $\mathcal{A}^*$, and we provide some first properties of this class, to one- or multi-dimensions.

The lower Karamata index of a one-dimensional distribution $G$ is defined as
\beam \label{eq.KPT.c.1}
K_G^- = \lim_{v \downarrow 1} \dfrac{-\log \overline{G^*}(v)}{\log v}\,. 
\eeam
The Karamata indexes were studied thoroughly in \cite[Subsestion 2.1]{bingham:goldie:teugels:1987}. As was mentioned also in \cite[p. 88]{cline:samorodnitsky:1994}, for any distribution $G$, with $\bG(x) > 0$ for any $x\in \bbr$, we have that $0\leq K_G^- \leq \infty$. Additionally, from \eqref{eq.KPT.c.1}, it is easy to see that if $G \in \mathcal{R}_{-\alpha}$, with $\alpha \in (0,\,\infty)$, then we obtain $K_G^- = \alpha$, while for any distribution $G$, with unbounded from above support, it holds $K_G^- \leq J_G^-$. Finally, we can easily check that if two distributions $G_1$ and $G_2$, are such that $\bG_1(x) \sim c\,\bG_2(x)$, for some $c \in (0,\,\infty)$, then we have $K_{G_1}^-=K_{G_2}^-$, while if $\Theta \sim G$, and $\lambda \in (0,\,\infty)$, with $\lambda\,\Theta \sim G'$, then we find $K_G^-=K_{G'}^-$.

In this way, in \cite{tang:yuan:2016} was introduced the class $\mathcal{A}^*$ as follows: We say that a one-dimensional distribution $G$ belongs to $\mathcal{A}^*$, if $G \in \mathcal{S}$ and $K_G^- >0$.

From the previous definition follow some interesting comments. At first, we notice that the most distributions $G$ with infinite right endpoint, are such that it holds $K_G^- =J_G^-$. Further, it is easy to see, due to inequalities in lower Karamata and Matuszewska indexes, that $\mathcal{A}^* \subsetneq \mathcal{A}$. However, the class $\mathcal{A} \setminus \mathcal{A}^*$, namely the class of subexponential distributions with $K_G^- =0$ and $J_G^->0$, is very poor. This, in combination with the fact that class $\mathcal{A}$ contains all the useful subexponential distributions, see relative discussions in \cite[Section 2]{tang:2006}, makes class $\mathcal{A}^*$ practically negligible smaller than class $\mathcal{S}$.

Now, we proceed to introduction of the multivariate analogue of $\mathcal{A}^*$, in similar way with the classes from Subsection 2.1.

\bde \label{def.KPT.c.1}
Let ${\bf X} \sim F$ be non-negative random vector, and let $A \in \mathscr{R}$ be some fixed set. We say that $F \in \mathcal{A}_A^*$, if $F_A \in \mathcal{A}^*$. Further, we denote
\beao
\mathcal{A}_\mathscr{R}^* = \bigcap_{A \in \mathscr{R}}\mathcal{A}_A^*\,.
\eeao
\ede

\bre \label{rem.KPT.c.1}
The class $\mathcal{A}_\mathscr{R}^*$, as also the rest of the classes in \eqref{eq.KLP.2.10}, has the following desired property: If ${\bf X} \sim F$, with $F \in \mathcal{A}_\mathscr{R}^*$, then all the non-negative and non-degenerate to zero linear combinations of components of ${\bf X}$, say $\sum_{j=1}^d l_j\,X_j$, have distribution from class $\mathcal{A}^*$. This can be verified, by considering as $A$, the set $A_2$ in Remark \ref{rem.KPT.3.1} below, from where is directly implied that, if $F \in \mathcal{A}_{A_2}^*$, then the $\sum_{j=1}^d l_j\,X_j$ follows distribution from class $\mathcal{A}^*$, for the $l_1,\,\ldots,\,l_d$ of set $A_2$. Therefore, from this and the definition of $\mathcal{A}_\mathscr{R}^* $, we find the desired result.
\ere  

The following proposition shows that class  $\mathcal{A}_\mathscr{R}^* $ contains class $MRV$.

\bpr \label{pr.KPT.2.1}
For any $\alpha \in (0,\,\infty)$, it holds $MRV(\alpha,\,\mu) \subsetneq \mathcal{A}_\mathscr{R}^* $.
\epr

\pr~
In proof of \cite[Proposition 4.14]{samorodnitsky:sun:2016} was shown that for any $A \in \mathscr{R}$, it holds $\mu(\partial A) =0$, and further $\mu(A) \in (0,\,\infty)$. Hence, for any $A \in \mathscr{R}$, since ${\bf X} \sim F \in MRV(\alpha,\,\mu)$, we have
\beam \label{eq.KPT.c.2}
\bF_A(x) = \PP\left({\bf X} \in x\,A \right) \sim \mu(A)\,\bG(x)\,, 
\eeam
with $G \in \mathcal{R}_{-\alpha}$, $\alpha \in (0,\,\infty)$. From \eqref{eq.KPT.c.2}, and the fact that $\mathcal{R}_{-\alpha}$ is closed 
with respect to strong-tail equivalence, see for example in \cite[Proposition 3.3(i)]{leipus:siaulys:konstantinides:2023}, we get $F_A \in \mathcal{R}_{-\alpha}$, 
and consequently $F_A \in \mathcal{A}^*$ (recall that since $F_A \in \mathcal{R}_{-\alpha}$, it follows that $K_{F_A}^- = \alpha >0$). Hence, $F \in \mathcal{A}_A^*$, 
and from the arbitrariness in choice of $A \in \mathscr{R}$, we obtain that $F \in \mathcal{A}_\mathscr{R}^*$.
~\halmos

\bre \label{rem.KPT.c.2}
From Proposition \ref{pr.KPT.2.1}, we can complete relation \eqref{eq.KLP.2.10} in the following way 
\beam \label{eq.KPT.c.3}
MRV\subsetneq (\mathcal{D} \cap \mathcal{A}^*)_\mathscr{R} \subsetneq \mathcal{A}^*_\mathscr{R} \subsetneq  \mathcal{S}_\mathscr{R}\subsetneq  \mathcal{L}_\mathscr{R}\,.
\eeam
Although $(\mathcal{D} \cap \mathcal{A})_\mathscr{R}$ is NOT proper subset of class $ \mathcal{A}_\mathscr{R}^*$, the $\mathcal{A}_\mathscr{R}^*$ contains almost all important distributions of $(\mathcal{D} \cap \mathcal{A})_\mathscr{R}$.
\ere

As we mentioned before, class $\mathcal{A}_\mathscr{R}^*$ is sufficiently richer than $MRV$. We shall discuss later and provide an example, which indicates this fact. However, for this we need the following proposition, that has its one merit, since gives a sufficient closure property condition of class $\mathcal{A}^*$ with respect to convolution. 

\bpr \label{pr.KPT.c.2.2}
Let $F_1,\,F_2$ be one-dimensional distributions with support on $\bbr$. If $F_1,\,F_2 \in \mathcal{A}^*$ and 
\beam \label{eq.KPT.c.4}
\overline{F_1*F_2}(x) \sim \bF_1(x) + \bF_2(x)\,.
\eeam
then it holds $F_1*F_2 \in \mathcal{A}^* $.
\epr  

\pr~
Firstly, since $F_1,\,F_2 \in \mathcal{A}^* \subsetneq \mathcal{A}$, and relation \eqref{eq.KPT.c.4} is true, by \cite[Corollary 3.1]{konstantinides:passalidis:2025b} is implied that $F_1 * F_2 \in \mathcal{A}$. So, it remains to show that $K_{F_1*F_2}^- >0$. For any $v> 1$ we have 
\beao
\overline{(F_1*F_2)^*}(v) = \limsup \dfrac{\overline{F_1*F_2}(v x) }{\overline{F_1*F_2}(x)} = \limsup \dfrac{\bF_1(v x)+ \bF_2(v x) }{\bF_1(x)+\bF_2(x)} \leq \max \left\{\overline{F_1^*}(v)\,,\; \overline{F_2^*}(v)\right\},
\eeao
where at the second step we used \eqref{eq.KPT.c.4}. Hence
\beao
K_{F_1*F_2}^- = \lim_{v\downarrow 1} \dfrac{-\log\left[\overline{(F_1*F_2)^*}(v)\right]}{\log v} \geq  \lim_{v\downarrow 1} \dfrac{-\log\left[\max \left\{\overline{F_1^*}(v)\,,\; \overline{F_2^*}(v)\right\}\right]}{\log v} >0\,,
\eeao
where the last step follows from the fact that $F_1,\,F_2 \in \mathcal{A}^*$.
~\halmos

In next example, we provide sufficient conditions with respect to marginal distributions, and their dependence structures, such that a multivariate distribution to belong to $\mathcal{A}_{A_3}^*$, with
\beam \label{eq.KPT.c.5}
A_3 = \left\{ {\bf y}\;:\;\dfrac 1d \sum_{i=1}^d l_i\,y_i > 1 \right\}\,,
\eeam
with $l_1,\,\ldots,\,l_d > 0$ and $l_1+ \cdots + l_d =1$.

\bexam \label{exam.KPT.c.2.1}
Let ${\bf X} \sim F$ be a non-negative random vector, with components $X_j \sim F_j \in \mathcal{A}^*$, for any $j = 1,\,\ldots,\,d$, and $\bF_j(x) \asymp \bF_n(x)$, for any $1\leq j \neq n \leq d$. Let us also assume that the $X_1,\,\ldots,\,X_d$, satisfy the following dependence structure, that was introduced by \cite{ko:tang:2008}, and contains the independence as special sub-case: There exists some $x_0 = x_0(d) > 0$ and some constant $C = C(d) \geq 0$, such that for any $x \geq x_0$ and any $2\leq n \leq d$ it holds 
\beam \label{eq.KPT.c.6}
\dfrac{\PP\left(\sum_{j=1}^{n-1} X_j > x-t\;\big|\;X_n = t \right)}{\PP\left(\sum_{j=1}^{n-1} X_j > x-t \right)} \leq M\,,
\eeam 
uniformly for $t \in [x_0,\,x]$. Next, we show that $F \in \mathcal{A}_{A_3}^*$. 

For the set $A_3$ given by relation \eqref{eq.KPT.c.5}, we have $I_{A_3}=\left(\dfrac 1{d},\,\ldots,\, \dfrac 1{d}\right)$. Hence, we obtain
\beao
\bF_{A_3}(x) = \PP\left({\bf X} \in x\,A \right) =\PP\left(\dfrac 1d \,\sum_{j=1}^{d} X_j > x \right)\,,
\eeao
for any $x>0$. Thus, from the conditions on the marginal distributions, and the dependence structure of \eqref{eq.KPT.c.6}, applying \cite[Theorem 3.1]{geng:liu:wang:2024}, we find that
\beam \label{eq.KPT.c.7}
\bF_{A_3}(x) = \PP\left(\sum_{j=1}^{d} \dfrac 1d\,X_j > x \right) \sim \sum_{j=1}^{d}\PP\left(\dfrac 1d\, X_j > x \right) \,,
\eeam
Let consider $X_j/d \sim F_j'$. So, since relation \eqref{eq.KPT.c.6} contains the independence as special sub-case, and due to $ F_j' \in \mathcal{A}^* $ and $\bF'_j(x) \asymp \bF'_n(x)$, for any $1 \leq j \neq n \leq d$, from \cite[Theorem 3.1]{geng:liu:wang:2024} follows that
\beam \label{eq.KPT.c.8}
\overline{F'_1 *F'_2*\,\cdots\,F'_d}(x) \sim \sum_{j=1}^{d}\PP\left( \dfrac{1}d\,X_j> x \right) = \sum_{j=1}^{d} \bF'_j(x) \,.
\eeam
From \eqref{eq.KPT.c.8}, and applying Proposition \ref{pr.KPT.c.2.2} inductively, we conclude that $F'_1 *F'_2*\,\cdots\,F'_d  \in \mathcal{A}^*$, and further by \eqref{eq.KPT.c.7} and \eqref{eq.KPT.c.8}, is implied that $F_{A_3} \in \mathcal{A}^*$.  Hence, we get $F \in \mathcal{A}_{A_3}^*$.  ~\halmos
\eexam

\bre \label{rem.KPT.c.3}
In the previous example we consider a weak dependence structure, and the condition of 
weak equivalence among the marginal distribution tails. In fact class $\mathcal{A}^*_\mathscr{R}$ or class $\mathcal{A}^*_A$, for some $A \in \mathscr{R}$, contains several cases, that either contain asymptotic dependence or the marginals have NOT asymptotically equivalent tails, or both. It is easy to check that if in \cite[Examples 4.1 -- 4.4]{konstantinides:liu:passalidis:2025} we additionally require 
for the marginals a positive Karamata index, then they belong to $\mathcal{A}^*_A$, 
for the set $A = A_2$ (see Remark \ref{rem.KPT.3.1} for $A_2$).

Furthermore, \cite[Example 4.5]{konstantinides:liu:passalidis:2025} belongs also to 
multivariate $\mathcal{A}^*_{A_2}$, for some set $A_2$, due to the fact that the rapidly varying distributions have infinite lower Karamata index. The last example does NOT contain any $MRV$ distribution.  
\ere

\subsection{Assumptions on the model}

In this subsection we provide two main assumptions for the discounted aggregate claims of relation \eqref{eq.KLP.1.2}. The first assumption simply says that the counting process $\{N(t)\,,\;t\geq 0\}$ is 
a renewal process, which is independent from all the other sources of randomness. 

\begin{assumption} \label{ass.KPT.2.1}
	We assume that the inter-arrival times $\{\theta_i:=\tau_i - \tau_{i-1}\,,\;i\in \bbn\}$, namely the interval times between arrivals of two successive claims, are i.i.d. non-negative random variables, which are independent from all the other sources of randomness. 
\end{assumption}

The following assumption relates to the dependence structure between the insurance
and financial risks. This dependence is inspired by \cite{asimit:jones:2008}, where 
was introduces a dependence structure through copulas, which later was generalized 
and used in several papers in continuous and discrete time, see for example in 
\cite{asimit:badescu:2010}, \cite{li:tang:wu:2010}, \cite{tang:yuan:2016}, \cite{yuan:lu:2023} among others. 

\begin{assumption} \label{ass.KPT.2.2}
	Let $A \in \mathscr{R}$, some fixed set. We assume that 
	\beao
	\{({\bf X}^{(i)}\,,\;e^{R(\tau_{i-1}) - R(\tau_i)})\,,\;i \in \bbn\}
	\eeao
	are i.i.d. copies of the general random vector $({\bf X}\,,\;e^{-R(\theta_1)})$, with 
	\beao
	{\bf X} \sim F\,, \qquad e^{-R(\theta_1)} \sim Q\,.
	\eeao
Further we assume that there exists some function $h\;:\;[0,\,\infty) \to (0,\,\infty)$, such that it holds
	\beao
	0 < \inf_{y \in E} h(y) \leq \sup_{y \in E} h(y) < \infty\,,
	\eeao
	where, $E$ represents a left neighborhood of the upper bound of $S(Q)$, such that it holds 
	\beam \label{eq.KPT.2.8}
	\PP\left( {\bf X} \in x\,A\;|\;e^{-R(s_1)} = y\,,\;\theta_1=s_1 \right) \sim h(y)\,\PP\left( {\bf X} \in x\,A \right) \,,
	\eeam 
	uniformly for any $y \in S(Q)$ and any $s_1 \in \Lambda$. 
\end{assumption}

It is worth to make some technical and intuitive comments on Assumption \ref{ass.KPT.2.2}.

\bre \label{rem.KPT.2.1}
The expression at the left side of relation \eqref{eq.KPT.2.8} is understood as follows
\beao
\lim_{\vep_1\wedge \vep_2 \to 0} \PP\left( {\bf X} \in x\,A\;|\;e^{-R(s_1)} \in [y-\vep_1,\,y+\vep_1]\,,\;\theta_1 \in [s_1-\vep_2,\,s_1+\vep_2] \right) \,.
\eeao
Furthermore, if $s_1 \notin \Lambda$, or if $s_1 \in \Lambda$ but $\PP\left( e^{R(s_1)} =y\,,\;\theta_1 =s_1 \right) =0$, for any $y \in \Delta$, for some open and non-empty set $\Delta$, then the conditional probability at the left side of relation \eqref{eq.KPT.2.8} is understood as unconditional, and hence $h(y)=1$, holds for any $y \in \Delta$, or any $y \in \bbr_+$ if $s_1 \notin \Lambda$.

The uniformity at relation \eqref{eq.KPT.2.8} means that it holds
\beao
\lim \sup_{s_1 \in \Lambda} \sup_{y \in S(Q)} \left|\dfrac{\PP\left( {\bf X} \in x\,A\;|\;e^{-R(s_1)} = y_1\,,\;\theta_1=s_1 \right) }{h(y)\,\PP\left( {\bf X} \in x\,A \right)} -1 \right|=0\,.
\eeao
\ere

\bre \label{rem.KPT.2.2}
In fact relation \eqref{eq.KPT.2.8} permits not only dependence between financial and 
insurance risks, but simultaneously dependence between the insurance risks, through the 
claim vectors, and the inter-arrival times. However, in this paper we consider only the 
independence assumption between the claims and the inter-arrival times, through Assumption \ref{ass.KPT.2.1}, and we do NOT extend our attention to time-dependent risk models. For more details about time-dependent risk models but with independent actuarial and financial risks, we refer to \cite{li:2012}, \cite{li:2016} among others. Hence, from Assumption \ref{ass.KPT.2.2} in combination with Assumption \ref{ass.KPT.2.1}, relation \eqref{eq.KPT.2.8} indicates that it holds 
\beao
\PP\left( {\bf X}^{(i)} \in x\,A\;|\;e^{R(s_{i-1})-R(s_i)} = y\right) \sim h(y)\,\PP\left( {\bf X} \in x\,A \right)\,,
\eeao 
for any $i \geq 1$, uniformly for any $y \in S(Q)$ and for any $s_{i-1} \in \Lambda$, 
with $s_{i-1} \leq s_i$. Further, by the fact that $\{({\bf X}^{(i)}\,,\;e^{R(\tau_{i-1}) - R(\tau_i)})\,,\;i \in \bbn\}$ are i.i.d., we can see that for any $i \in \bbn$ the insurance risk ${\bf X}^{(i)}$ is dependent on its financial returns $\{R(t)\,,\;t\geq 0\}$ only over the interval $(\tau_{i-1}\,,\;\tau_i]$, while the ${\bf X}^{(i)}$ is independent of the increments of the L\'{e}vy process over the intervals $[0\,,\;\tau_{i-1}]$ and $(\tau_i\,,\;\infty]$, which seems reasonable in combination with the assumption that the logarithmic returns are described by a  L\'{e}vy process. 
\ere

\bre \label{rem.KPT.2.3}
We can observe that the dependence in relation \eqref{eq.KPT.2.8} contains as special cases either the independence between ${\bf X}$ and $R(t)$, or the independence between ${\bf X}$ and $\theta$, or the independence among ${\bf X}$, $R(t)$ and $\theta$. Further, relation \eqref{eq.KPT.2.8} indicates the existence of a weak dependence structure between the insurance and financial risks, in the sense of asymptotic independence, see in \cite{asimit:furman:tang:vernic:2011} for related discussion about the static structure of this dependence. From the static structure of this dependence, in combination with Assumption \ref{ass.KPT.2.1}, we can find out that the dependence in relation \eqref{eq.KPT.2.8} is rich enough, see \cite[Sec. 3]{li:tang:wu:2010} and \cite[Sec. 2.1]{cui:wang:2025} for examples of static form  of this dependence. Therefore, from the one hand side, we obtain a general enough structure of weak dependence between the two risks, that contains the independence as special case, and from the other hand side, since in Section \ref{sec.KLP.3} we consider that the L\'{e}vy processes are non-negative, namely insurer invest only on risk-free assets, we expect that the dependencies in practice will be rather weak than strong.   
\ere

\bre \label{rem.KPT.2.4}
Finally, we show that the insurance and financial risks belong to $\mathcal{A}_A$, when $F \in \mathcal{A}_A$, for some $A \in \mathscr{R}$, $\{R(t)\,,\;t\geq 0\}$ is a non-negative L\'{e}vy process, and the Assumptions  \ref{ass.KPT.2.1}, \ref{ass.KPT.2.2} hold. Let $\{R^*(t)\,,\;t\geq 0\}$ be a new c\'{a}dl\'{a}g process, that is independent from the other sources of randomness and the distribution of $e^{-R^*(t)}$ is given by
\beam \label{eq.KPT.2.10} 
\PP\left( e^{-R^*(t)} \in dy\right)= h(y)\,\,\PP\left( e^{-R(t)} \in dy \right)\,,
\eeam 
for any $t\geq 0$. It is easy to check that from relation \eqref{eq.KPT.2.8}, after integration with respect to $y$ and $s_1$ in $S(Q)$ and $\Lambda$ that 
\beao
\E\left[h\left(e^{-R(\theta_1)} \right) \right]=1\,.
\eeao 
This, in combination with relation \eqref{eq.KPT.2.10}, implies that the $e^{-R^*(\theta_1)}$ has a proper distribution. Hence we obtain  
\beam \label{eq.KPT.2.11} \notag
&&\PP\left( {\bf X}\,e^{-R(\theta_1)} \in x\,A\right)= \int_{\Lambda}\,\PP\left( {\bf X}\,e^{-R(s_1)} \in x\,A \right)\,\PP(\theta_1 \in ds_1)\\[2mm] \notag
&&= \int_{\Lambda} \int_0^1\,\PP\left( {\bf X} \in \dfrac xy\,A \;\Big|\;e^{-R(s_1)} =y\right)\,\PP\left( e^{-R(s_1)} \in dy \right)\,\PP(\theta_1 \in ds_1)\\[2mm] \notag
&& \sim \int_{\Lambda} \int_0^1 h(y)\,\PP\left( {\bf X} \in \dfrac xy\,A \right)\,\PP\left( e^{-R(s_1)} \in dy \right)\,\PP(\theta_1 \in ds_1)\\[2mm] \notag
&&=\int_{\Lambda} \int_0^1 \PP\left( {\bf X} \in \dfrac xy\,A \right)\,\PP\left( e^{-R^*(s_1)} \in dy \right)\,\PP(\theta_1 \in ds_1)\\[2mm] 
&&=\int_{\Lambda} \PP\left( {\bf X}\,e^{-R^*(s_1)} \in  x\,A \right)\,\PP(\theta_1 \in ds_1)=\PP\left( {\bf X}\,e^{-R^*(\theta_1)} \in  x\,A \right)\,,
\eeam 
where at the first step we used Assumption \ref{ass.KPT.2.1}, at the third step we took into account relation \eqref{eq.KPT.2.8}, and at fourth step we employ relation  relation \eqref{eq.KPT.2.10}. But, since $e^{-R^*(\theta_1)}$ belongs to the interval $[0,\,1]$ and ${\bf X} \sim F \in \mathcal{A}_A$, with ${\bf X} $ and $e^{-R^*(\theta_1)}$ independent each other, from \cite[Theorem 4.1]{konstantinides:passalidis:2024h} we find ${\bf X}\,e^{-R^*(\theta_1)} \in \mathcal{A}_A$. Hence, from relation \eqref{eq.KPT.2.11}, in combination with \cite[Proposition 3.2]{konstantinides:passalidis:2024h} we have ${\bf X}\,e^{-R(\theta_1)} \in \mathcal{A}_A$.
\ere

\section{Weakly Dependent Insurance and Financial Risks}\label{sec.KLP.3}

\subsection{Main result} \label{sec.KLP.3.3.1}

Here we give the main result for the asymptotic behavior of the entrance probability 
of the discount claims ${\bf D}(\infty)$ into some rare set $x\,A$, with $A \in \mathscr{R}$, when the insurer makes risk-free investments, and exists the weak dependence structure of Assumption \ref{ass.KPT.2.2} between the insurance and financial risks. We also notice that we assume that the claim vector distribution $F$, belongs to class $\mathcal{A}_A^*$.  At first we present our result and then we begin some discussion about the generality of the assumptions and its applicability in actuarial practice. 

\bth \label{th.KPT.3.1}
Let  $A \in \mathscr{R}$ be some fixed set and let consider the discounted aggregate claims from relation \eqref{eq.KLP.1.3}. We suppose that the Assumptions \ref{ass.KPT.2.1} and \ref{ass.KPT.2.2} are true and additionally the process $\{R(t)\,,\;t\geq 0\}$ is non-negative and non-degenerate to zero, L\'{e}vy process. We also assume that $F \in \mathcal{A}_A^*$. Then it holds
\beam \label{eq.KPT.3.1} 
\PP\left( {\bf D}(\infty) \in x\,A\right) \sim \int_{0}^{\infty} \PP\left( {\bf X} \,e^{-R(s)} \in x\,A \right)\,\lambda(ds)\,.
\eeam
\ethe

\bre \label{rem.KPT.3.1}
Relation \eqref{eq.KPT.3.1} was established under various assumptions for the dependencies and the distribution classes either for the claim vectors, or for the financial risks, see for example in \cite{chen:konstantinides:passalidis:2025}, 
\cite{konstantinides:passalidis:2024j}. In Theorem \ref{th.KPT.3.1}, in spite of the restriction to non-negative L\'{e}vy processes, and the dependence Assumption \ref{ass.KPT.2.2} between the two risks, we get the product distribution of ${\bf X}\,e^{-R(\theta_1)}$ to belong to $\mathcal{A}_A$, as we see in Remark \ref{rem.KPT.2.4}. 

In general, the expression of relation \eqref{eq.KPT.3.1} discloses the multivariate linear single big jump principle of the discounted aggregate claims. In one-dimensional case, the sets of the form $A=(b,\,\infty)$, with $b>0$, belong to family $\mathscr{R}$. Hence, the set $A=(1,\,\infty)$ in relation \eqref{eq.KPT.3.1} shows a classical asymptotic estimation for the one-dimensional behavior of the distribution tail of the discounted aggregate claims when the claim distribution is subexponential (or some subclass). Even in such special one-dimensional case, Theorem \ref{th.KPT.3.1} produces  new result. 

Further more, two other sets from the family  $\mathscr{R}$, with huge practical interest, are given below 
\beao
A_1:=\{{\bf y}\;:\; y_j > b_j\,,\;\exists \;j=1,\,\ldots,\,d\}\,,
\eeao
for $b_1,\,\ldots,\,b_d>0$ and 
\beao
A_2:=\left\{{\bf y}\;:\; \sum_{j=1}^d l_j\,y_j > b \right\}\,,
\eeao
with $l_1,\,\ldots,\,l_d\geq 0$, $b>0$ and
\beao
\sum_{j=1}^d l_j =1\,.
\eeao                                                      
\ere

\bre \label{rem.KPT.3.2}
If in  Theorem \ref{th.KPT.3.1} we consider that $F \in \mathcal{A}_{\mathscr{R}}^*$, and Assumption \ref{ass.KPT.2.2} holds for any  $A \in \mathscr{R}$, then relation \eqref{eq.KPT.3.1} also holds for any  $A \in \mathscr{R}$. Further, relation \eqref{eq.KPT.3.1} seems to be difficult for computation, because of the generality of assumptions. However,  when the set $A$ and the components ${\bf X}$ dependence structures are determined, then relation \eqref{eq.KPT.3.1} can be reduced into more direct asymptotic expressions, which are easily computable.
\ere
Furthermore, as we see later in Corollary \ref{cor.KLP.5.1} (under weaker conditions on L\'{e}vy process $\{R(t)\,,\;t\geq 0\}$), when the requirement $F \in \mathcal{A}_A^*$ is replaced by $F \in MRV$, this implies an even more immediate asymptotic expression in comparison with \eqref{eq.KPT.3.1}.

\subsection{Proof of Theorem \ref{th.KPT.3.1} }\label{sec.KLP.3.3.2}

Before the proof of theorem, we need some preliminary lemmas. At first we provide the following 'static' form of the Assumption \ref{ass.KPT.2.2}, which is necessary to the preliminary lemmas.

\begin{assumption} \label{ass.KPT.3.1}
Let $A \in \mathscr{R}$ be a fixed set. Let $({\bf X}^{(1)}\,,\;W_1),\,\ldots,\,({\bf X}^{(n)}\,,\;W_n)$, with $n \in \bbn$, be a sequence of non-negative, independent, but not necessarily identically distributed, random pairs. We suppose that for any $i=1,\,\ldots,\,n$, it holds ${\bf X}^{(i)}\sim F_i$ and there exists a constant $b_i \in (0,\,\infty)$, such that $\PP(0\leq W_i \leq b_i)=1$, $\PP(W_i=0)<1$ hold and there exists a function $h_i\;:\;[0,\,\infty) \to (0,\,\infty)$, such that the inequalities $0< \inf_{y \in E_i} h_i(y) \leq \sup_{y \in E_i} h_i(y) < \infty$ hold, where $E_i$ is a left neighborhood of $b_i$, for which we obtain
	\beam \label{eq.KPT.3.3} 
	\PP\left( {\bf X}^{(i)} \in x\,A\;|\;W_i=y \right) \sim h_i(y)\,\PP\left( {\bf X}^{(i)}\in x\,A \right)\,,
	\eeam
	uniformly for $y \in S(W_i)$. 
\end{assumption}

Let us notice that in Assumption \ref{ass.KPT.3.1} the uniformity of \eqref{eq.KPT.3.3} is understood in the sense of
\beao
\lim \sup_{y \in S(W_i)} \left|\dfrac{\PP\left( {\bf X}^{(i)} \in x\,A\;|\;W_i=y \right)}{h_i(y)\,\PP\left( {\bf X}^{(i)}\in x\,A \right)} -1 \right|=0\,.
\eeao

In the following lemma, when we say that the pair $({\bf X}\,,\;W)$ satisfies Assumption \ref{ass.KPT.3.1}, we have in mind that ${\bf X} \sim F$, $\PP(0\leq W \leq b)=1$, for some $b \in (0,\,\infty)$, $\PP(W=0)<1$ and relation \eqref{eq.KPT.3.3} is true with a function $h$, that is bounded away from zero and infinity over some left area $E$ of $b$. The following lemma represents a multivariate extension of \cite[Lem. 4.3]{tang:yuan:2016}.

\ble \label{lem.KPT.3.1}
Let $A \in \mathscr{R}$ be some fixed set. We assume that $({\bf X}\,,\;W)$ is a non-negative random pair, that satisfies Assumption \ref{ass.KPT.3.1}, with ${\bf X} \sim F \in \mathcal{S}_A$. Let us consider another non-negative random vector ${\bf Z} \sim V \in \mathcal{S}_A $, which is independent of $({\bf X}\,,\;W)$. If it holds either $F(x\,A)= O[V(x\,A)]$, or $V(x\,A) = O[F(x\,A)]$, then we obtain
\beam \label{eq.KPT.3.4} 
\PP\left( [{\bf X}+{\bf Z} ]\,W \in x\,A \right) \sim \PP\left( {\bf X}\,W \in x\,A \right) + \PP\left( {\bf Z}\,W \in x\,A \right)\,.
\eeam
\ele   

\pr~
Let $Z_A$ the random variable of relation \eqref{eq.KPT.2.2} following distribution $V_A$ and similarly remember that
\beam \label{eq.KLP.3.1b}
X_A := \{ u\;:\; {\bf X} \in u\,A \} \sim F_A\,.
\eeam
Hence, $F_A \in \mathcal{S}$, $V_A \in \mathcal{S}$ and either 
\beao
\bF_A(x)=O[\bV_A(x)]\,,
\eeao 
or  
\beao
\bV_A(x)=O[\bF_A(x)]
\eeao 
hold, then by \cite[Prop. 2.4]{konstantinides:passalidis:2024g} 
and further by \cite[Lem. 4.3]{tang:yuan:2016} we obtain
\beao
&&\PP\left( [{\bf X}+{\bf Z} ]\,W \in x\,A \right) \leq \PP\left( X_A\,W + Z_A\,W > x \right)\\[2mm]
&& \sim \PP\left( X_A\,W > x \right) + \PP\left(Z_A\,W > x \right)=\PP\left({\bf X}\,W \in x\,A \right)+ \PP\left( {\bf Z}\,W \in x\,A \right)\,,
\eeao 
that gives the desired upper bound of relation \eqref{eq.KPT.3.4}.

So, we estimate now the lower bound of relation \eqref{eq.KPT.3.4}. From the fact that the set $A$ is increasing, and the ${\bf X}$, $W$, {\bf Z} are non-negative, employing the Bonferroni inequality, we obtain
\beam \label{eq.KPT.3.5} 
&&\PP\left( [{\bf X}+{\bf Z} ]\,W \in x\,A \right) = \int_0^b \PP\left( {\bf X}+ {\bf Z}  \in \dfrac xy\,A\;|\;W=y \right)\,\PP\left(W \in dy \right)\\[2mm] \notag
&& \geq \int_0^b \PP\left( {\bf X} \in \dfrac xy\,A\;|\;W=y \right)\,\PP\left(W \in dy \right)+ \int_0^b \PP\left({\bf Z} \in \dfrac xy\,A\right)\,\PP\left(W \in dy \right)\\[2mm] \notag
&& - \int_0^b \PP\left( {\bf X}  \in \dfrac xy\,A\,,\; {\bf Z} \in \dfrac xy\,A\,,\;|\;W=y \right)\,\PP\left(W \in dy \right) \\[2mm] \notag
&&\sim \int_0^b h(y)\,\PP\left( {\bf X} \in \dfrac xy\,A \right)\,\PP\left(W \in dy \right)+ \PP\left({\bf Z}\,W \in x\,A\right)\\[2mm] \notag
&& - \int_0^b \PP\left( {\bf Z}  \in \dfrac xy\,A \right)\,h(y)\,\PP\left( {\bf X} \in \dfrac xy\,A \right)\,\PP\left(W \in dy \right) \,,
\eeam
where at the last step in relation \eqref{eq.KPT.3.5} we used Assumption \ref{ass.KPT.3.1}, in combination with the fact that ${\bf Z}$ is independent of $({\bf X},\,W)$. 

Let consider now a new random variable $W_h$, that is independent of all the sources 
of randomness with distribution
\beam \label{eq.KPT.3.6} 
\PP\left( W_h \in dy \right) = h(y)\,\PP\left(W \in dy \right)\,,
\eeam
for any $y \in S(W)$, which is a proper distribution, since $\E[h(W)] = 1$, that follows immediately from Assumption \ref{ass.KPT.3.1}, through integration on both sides of relation \eqref{eq.KPT.3.3} with respect to $\PP\left(W \in dy \right)$ over the whole $S(W)$. Hence, by relations \eqref{eq.KPT.3.5} and \eqref{eq.KPT.3.6} we find that
\beam \label{eq.KPT.3.7} \notag
\PP\left( [{\bf X}+{\bf Z} ]\,W \in x\,A \right) &\gtrsim& \PP\left( {\bf X}\,W_h \in x\,A \right) + \PP\left( {\bf Z}\,W \in x\,A \right)\\[2mm]
&-& \PP\left( {\bf Z}\,W_h \in x\,A \right)\,\PP\left( {\bf X}\,W_h \in x\,A \right)\,.
\eeam
From the fact that $W$ represent a bounded from above random variable, we can use 
\cite[Lemma 3.1]{chen:xu:cheng:2019} to obtain
\beam  \label{eq.KPT.3.8} 
\PP\left( {\bf X}\,W_h \in x\,A \right) \sim \PP\left( {\bf X}\,W \in x\,A \right)\,,
\eeam
while from the other hand side we find that it holds
\beam \label{eq.KPT.3.9} 
\PP\left( {\bf Z}\,W_h \in x\,A \right) \asymp \PP\left( {\bf Z}\,W \in x\,A \right)\,.
\eeam
Indeed, relation \eqref{eq.KPT.3.9} is true, since for some left area $E$ of $b$, the 
function $h$ is bounded away from zero and infinity on some neighborhood $E$, hence by relation  \eqref{eq.KPT.3.6} we obtain
\beao
&&\PP\left( {\bf Z}\,W_h \in x\,A \right) = \int_0^b h(y)\,\PP\left({\bf Z} \in \dfrac xy\,A \right)\,\PP\left(W \in dy \right)\\[2mm]
&&\asymp \int_E h(y)\,\PP\left({\bf Z} \in \dfrac xy\,A \right)\,\PP\left(W \in dy \right) \asymp \int_E \PP\left({\bf Z} \in \dfrac xy\,A \right)\,\PP\left(W \in dy \right)\\[2mm]
&&\asymp \int_0^b \PP\left({\bf Z} \in \dfrac xy\,A \right)\,\PP\left(W \in dy \right) = \PP\left( {\bf Z}\,W \in x\,A \right)\,,
\eeao
where the second and fourth step are easily verified, by separating the interval $[0,\,b]$ in $E$ and $E^c$. Therefore, from relations \eqref{eq.KPT.3.7}, \eqref{eq.KPT.3.8} and \eqref{eq.KPT.3.9} we have
\beao
\PP\left( [{\bf X}+{\bf Z} ]\,W \in x\,A \right) \gtrsim \PP\left( {\bf X}\,W \in x\,A \right) + \PP\left( {\bf Z}\,W \in x\,A \right)\,, 
\eeao 
which gives the desired lower bound of relation \eqref{eq.KPT.3.4} and the proof is completed. 
~\halmos

The following lemma refers to multivariate linear single big jump principle of scale mixture sums, and its first part represent a reformulation of the \cite[Lem. 5.1(i)]{konstantinides:passalidis:2024h}, while at the second step we show the multivariate  subexponentiality for this scale mixture sum.

\ble \label{lem.KPT.3.2}
Let $A \in \mathscr{R}$ be some fixed set. We consider the sequence ${\bf X}^{(1)},\,\ldots,\,{\bf X}^{(n)}$ with $n \in \bbn$, of non-negative, independent random vectors, with distributions $F_1,\,\ldots,\,F_n \in \mathcal{L}_A$ respectively, and we suppose that $F_i(x\,A) \asymp F(x\,A)$, for any $i=1,\,\ldots,\,n$, with $F \in \mathcal{S}_A$. Let us assume that the $\Theta_1,\,\ldots,\,\Theta_n$ are arbitrarily dependent, non-negative and non-degenerate to zero, bounded from above random variables, which are independent of  ${\bf X}^{(1)},\,\ldots,\,{\bf X}^{(n)}$. Then we obtain:
\begin{enumerate}
	\item[(i)] 
	\beam \label{eq.KPT.3.10}
	\PP\left( \sum_{i=1}^n \Theta_i\,{\bf X}^{(i)} \in x\,A \right) \sim \sum_{i=1}^n \PP\left( \Theta_i\,{\bf X}^{(i)} \in x\,A \right)\,.
		\eeam
	\item[(ii)]
	The distribution of $\sum_{i=1}^n \Theta_i\,{\bf X}^{(i)}$ belongs to $\mathcal{S}_A$.
\end{enumerate}
\ele 

\pr~
\begin{enumerate}
	\item[(ii)] 
In order to show that the distribution of $\sum_{i=1}^n \Theta_i\,{\bf X}^{(i)}$ belongs to $\mathcal{S}_A$, it is enough to show that the random variable
	\beao
	P_A:= \sup\left\{u\;:\; \sum_{i=1}^n \Theta_i\,{\bf X}^{(i)} \in u\,A \right\} =\sup_{{\bf p} \in I_A} {\bf p}^T\,\sum_{i=1}^n \theta_i\,{\bf X}^{(i)} \sim M_A
	\eeao
is such that it holds $M_A \in \mathcal{S}$. Let consider an independent copies of 
	\beao
	(\Theta_1,\,\ldots,\,\Theta_n,\,{\bf X}^{(1)},\,\ldots,\, {\bf X}^{(n)})\,,
	\eeao 
let denote by $(\Theta_1^*,\,\ldots,\,\Theta_n^*,\,{\bf X}^{(1)*},\,\ldots,\, {\bf X}^{(n)*})$. Then the 
	\beao
	P_A^*:= \sup\left\{u\;:\; \sum_{i=1}^n \Theta_i^*\,{\bf X}^{(i)*} \in u\,A \right\} \,, =\sup_{{\bf p} \in I_A} {\bf p}^T\,\sum_{i=1}^n \theta_i^*\,{\bf X}^{(i)*} 
	\eeao
has also distribution $M_A$ and is independent of $P_A$. Thus, it is enough to show that 
	$\PP(P_A +P_A^* >x) \sim 2\,\PP(P_A> x)$.
	
	For 
	\beao 
	X_A^{(i)}:= \sup \{u\;:\;{\bf X}^{(i)} \in u\,A\}\,,\; X_A^{(i)*}:= \sup \{u\;:\;{\bf X}^{(i)*} \in u\,A\}\,,
	\eeao
	we find 
	\beam \label{eq.KPT.3.13}
	&&\PP(P_A +P_A^* >x) = \PP\left( \sum_{i=1}^n \Theta_i\,X_A^{(i)} + \sum_{i=1}^n \Theta_i^*\, X_A^{(i)*} >x\right) \\[2mm] \notag
	&& \sim \PP\left( \sum_{i=1}^n \Theta_i\,X_A^{(i)} >x\right) + \PP\left( \sum_{i=1}^n \Theta_i^*\, X_A^{(i)*} >x\right) \\[2mm] \notag 
	&& \sim \sum_{i=1}^n \PP\left( \Theta_i\,X_A^{(i)} >x\right) + \sum_{i=1}^n\PP\left(  \Theta_i^*\, X_A^{(i)*} >x\right) =2\,\sum_{i=1}^n \PP\left( \Theta_i\,{\bf X}^{(i)} \in x\,A \right) \\[2mm] \notag 
	&& \sim 2\,\PP\left( \sum_{i=1}^n  \Theta_i\,{\bf X}^{(i)} \in x\,A\right)=2\,\PP\left( P_A>x\right)\,,
	\eeam
	where at the first step we used \cite[Prop. 2.4]{konstantinides:passalidis:2024g}, at the second step we used \cite[Lem. 4.4]{tang:yuan:2016}, namely that the distribution of 
	\beao
	\sum_{i=1}^n \Theta_i\,X_A^{(i)}
	\eeao
	belongs to $\mathcal{S}$, at the third step we employed \cite[Th. 1]{tang:yuan:2014}, while at the penultimate step we took into account relation \eqref{eq.KPT.3.10}. Hence, by \eqref{eq.KPT.3.13} follows that $P_A \sim M_A \in \mathcal{S}$. ~\halmos
\end{enumerate}

The following lemma it has also its own merit, since give the multivariate linear single big jump principle of the ${\bf S}_n$, which given in \eqref{eq.KPT.3.15}. We see that in the special case when $A=(1,\,\infty)$, this lemma  is reduced to \cite[Th. 3.5]{tang:yuan:2016}.

\ble \label{lem.KPT.3.3}
Let $A \in \mathscr{R}$ be some fixed set. We consider that the $({\bf X}^{(1)},\,W_1),\,\ldots,\,({\bf X}^{(n)},\,W_n)$ satisfy Assumption \ref{ass.KPT.3.1}, and further for any $i=1,\,\ldots,\,n$ there exists $c_i \in (0,\,\infty)$ such that it  holds 
\beao
F_i(x\,A) \sim c_i\,F(x\,A)\,,
\eeao 
with $F \in \mathcal{S}_A$. Then we obtain 
\beam \label{eq.KPT.3.14}
\PP\left( {\bf S}_n \in x\,A \right) \sim \sum_{i=1}^n \PP\left( {\bf X}^{(i)}\,\prod_{j=i}^n W_j \in x\,A\right)\,,
\eeam
where 
\beam \label{eq.KPT.3.15}
{\bf S}_n := \sum_{i=1}^n  {\bf X}^{(i)}\,\prod_{j=i}^n W_j \,.
\eeam
\ele

\pr~
At first we mention that the upper bound in relation \eqref{eq.KPT.3.14} can be verified through \cite[Prop. 2.4]{konstantinides:passalidis:2024g} and \cite[Th. 3.5]{tang:yuan:2016}, however here we follow a direct alternative approach for the proof of relation \eqref{eq.KPT.3.14}, that leads to the desired asymptotic estimation, without the verification of the upper and lower bounds separately. We observe also that it follows $F_i \in \mathcal{S}_A$, for any $i=1,\,\ldots,\,n$, by the closure property 
of $\mathcal{S}_A$ with respect to strong-tail equivalence, see in \cite[Prop. 4.12(a)]{samorodnitsky:sun:2016}.

We use induction technique to show relation \eqref{eq.KPT.3.14}. For $n=1$, relation \eqref{eq.KPT.3.14} holds trivially as equality.

Let assume that \eqref{eq.KPT.3.14} holds for some $k \in \bbn$. We show that it also holds for $k+1$. At first we define the $W_{h,i}$, for $i=1,\,\ldots,\,k$, that are independent each other, independent from all other sources of randomness random variables, with proper distributions 
\beao
\PP\left(W_{h,i} \in dy\right) = h_i(y)\,\PP\left( W_i \in dy\right)\,,
\eeao 
for any $y \in S(W_i)$. Hence, we find out that for any $i=1,\,\ldots,\,k$ it holds
\beam \label{eq.KPT.3.16}  \notag
&&\PP\left( {\bf X}^{(i)}\,\prod_{j=i}^k W_j \in x\,A\right)= \int_0^{b_i} \PP\left(  {\bf X}^{(i)}\,\prod_{j=i+1}^k W_j \in \dfrac xy \,A \;\Big|\;W_i = y \right)\,\PP\left( W_i \in dy\right) \\[2mm]  \notag
&& = \int_0^{b_i} \int_{0}^{b_{i+1}\cdot \ldots \cdot b_k}\PP\left(  {\bf X}^{(i)} \in \dfrac x{y\,m} \,A \;\Big|\;W_i = y \right)\,\PP\left( \prod_{j=i+1}^k W_j \in dm\right)\,\PP\left( W_i \in dy\right) \sim\\[2mm] 
&& \int_0^{b_i} \int_{0}^{b_{i+1}\cdot \ldots \cdot b_k} h_i(y)\,\PP\left( {\bf X}^{(i)} \in \dfrac x{y\,m} \,A \right)\,\PP\left( \prod_{j=i+1}^k W_j \in dm\right)\,\PP\left( W_i \in dy\right) \\[2mm] \notag
&&=  \int_0^{b_i} \PP\left(  {\bf X}^{(i)}\,\prod_{j=i+1}^k W_j \in \dfrac x{y} \,A \right)\,\PP\left( W_{h,i} \in dy\right)=\PP\left(  {\bf X}^{(i)}\,W_{h,i}\,\prod_{j=i+1}^k W_j \in x\,A \right)\,,
\eeam 
where at the third step we used relation \eqref{eq.KPT.3.3} through the dominated convergence theorem (remember the uniformity of \eqref{eq.KPT.3.3}). Hence, by the induction assumption, in combination with relation \eqref{eq.KPT.3.16} we find
\beao
\PP\left( {\bf S}_k \in x\,A \right) \sim \sum_{i=1}^k \PP\left( {\bf X}^{(i)}\,\prod_{j=i}^k W_j \in x\,A\right)  \sim \sum_{i=1}^k \PP\left( {\bf X}^{(i)}\,W_{h,i}\,\prod_{j=i+1}^k W_j \in x\,A\right)\,, 
\eeao  
so from the last relation, putting as random weights the 
\beao
\Theta_i= W_{h,i}\,\prod_{j=i+1}^k W_j\,,
\eeao 
that are bounded from above by $(b_i\cdot \ldots \cdot b_k)$ and independent of ${\bf X}^{(1)},\,\ldots,\,{\bf X}^{(k)}$, from Lemma \ref{lem.KPT.3.2}(i) we obtain
\beam \label{eq.KPT.3.17}  \notag
\PP\left( {\bf S}_k \in x\,A \right) \sim \sum_{i=1}^k \PP\left( {\bf X}^{(i)}\,W_{h,i}\,\prod_{j=i+1}^k W_j \in x\,A\right) \sim \PP\left( \sum_{i=1}^k {\bf X}^{(i)}\,W_{h,i}\,\prod_{j=i+1}^k W_j \in x\,A\right),\\
\eeam 
Since the last term of relation \eqref{eq.KPT.3.17} belongs to $\mathcal{S}_A$, due to Lemma  \ref{lem.KPT.3.2}(ii). Therefore, from closure property of $\mathcal{S}_A$ with respect to strong tail equivalence, we find that the distribution of ${\bf S}_k$ belongs to $\mathcal{S}_A$.

Now we examine separately two cases. In first case, we consider that 
\beao
\max_{1 \leq i \leq k} \prod_{j=1}^k b_j \leq 1\,,
\eeao 
and consequently, since the set $A$ is increasing set, we obtain 
\beam \label{eq.KPT.3.18} 
\PP\left( {\bf S}_k \in x\,A \right) \leq \PP\left(\sum_{i=1}^k  {\bf X}^{(i)} \in x\,A\right) \sim \sum_{i=1}^k \PP\left( {\bf X}^{(i)} \in x\,A\right)= O\left[\PP\left( {\bf X}^{(k+1)} \in x\,A\right) \right]\,,
\eeam
where at the second step we used Lemma \ref{lem.KPT.3.2}(i), with degenerate to unit weights. 

In second case, we consider that there exists some $1\leq l \leq k$, such that 
\beao
\prod_{j=l}^k b_j > 1\,,
\eeao
and consequently we find
\beao
&&\PP\left( {\bf S}_k \in x\,A \right) \geq \PP\left({\bf X}^{(l)}\,\prod_{j=l}^k W_j \in x\,A\,,\; \prod_{j=l}^k W_j > 1 \right) \\[2mm] 
&&\geq \PP\left({\bf X}^{(l)} \in x\,A\,,\; W_l\,\prod_{j=l+1}^k W_j > 1 \right) \\[2mm] 
&&= \int_0^{b_l} \PP\left({\bf X}^{(l)} \in x\,A\,,\;y\,\prod_{j=l+1}^k W_j >1 \,,\;|\;  W_l =y \right) \,\PP\left( W_l \in dy \right) \\[2mm] \notag
&&= \int_0^{b_l} \int_1^{b_{l+1}\ldots b_k} \PP\left({\bf X}^{(l)} \in x\,A\,,\;|\;  W_l =y \right)\,\PP\left(y\,\prod_{j=l+1}^k W_j \in dm \right) \,\PP\left( W_l \in dy \right) \\[2mm] \notag
&&\sim \int_0^{b_l} h_l(y) \PP\left({\bf X}^{(l)} \in x\,A \right)\,\PP\left(y\,\prod_{j=l+1}^k W_j >1 \right)\,\PP\left( W_l \in dy \right)  \\[2mm] 
&&= \PP\left({\bf X}^{(l)} \in x\,A\right)\,\PP\left(W_{h,l}\,\prod_{j=l+1}^k W_j >1 \right)\gtrsim C\,\PP\left( {\bf X}^{(k+1)} \in x\,A\ \right) \,,
\eeao
where at the last step the constant $C>0$ stems from the assumption that the $F_i$ are strong-tailed equivalent in combination with the fact that 
\beao
\prod_{j=l}^k b_j >1\,.
\eeao 
Thus, in this case it holds
\beam \label{eq.KPT.3.19} 
\PP\left( {\bf X}^{(k+1)} \in x\,A \right) =O\left[ \PP\left({\bf S}_k \in x\,A \right)\right]\,,
\eeam

Hence, from relations  \eqref{eq.KPT.3.18} and  \eqref{eq.KPT.3.19} and the fact that the ${\bf S}_k$ and ${\bf X}^{(k+1)}$ are independent with distributions from class $\mathcal{S}_A$, we can use Lemma Lemma  \ref{lem.KPT.3.1} at the second step of the following relation
\beam \label{eq.KPT.3.20} \notag
&&\PP\left( {\bf S}_{k+1} \in x\,A \right) = \PP\left(({\bf X}^{(k+1)} + {\bf S}_{k})\,W_{k+1} \in x\,A\right)  \\[2mm] \notag
&&\sim \PP\left( {\bf X}^{(k+1)}\,W_{k+1} \in x\,A\right)+ \PP\left( {\bf S}_{k}\,W_{k+1} \in x\,A\right)  \\[2mm] 
&&\sim \PP\left( {\bf X}^{(k+1)}\,W_{k+1} \in x\,A\right)+ \int_0^{b_{k+1}}\PP\left( {\bf S}_{k} \in \dfrac xy\,A\right)\,\PP\left( W_{k+1} \in dy \right) \\[2mm] \notag
&&= \PP\left( {\bf X}^{(k+1)}\,W_{k+1} \in x\,A\right)+ \sum_{i=1}^{k} \PP\left( {\bf X}^{(i)}\,\prod_{j=i}^{k+1} W_j \in x\,A\right)\\[2mm] \notag
&&= \sum_{i=1}^{k+1} \PP\left( {\bf X}^{(i)}\,\prod_{j=i}^{k+1} W_j \in x\,A\right) \,,
\eeam
where at the fourth step we used induction assumption, through the dominated convergence theorem. Relation \eqref{eq.KPT.3.20} completes the induction step. 
~\halmos

The following lemma is instrumental for the proof of Theorem \ref{th.KPT.3.1}. In fact this lemma says that under some stricter conditions in Lemma \ref{lem.KPT.3.3}, relation \eqref{eq.KPT.3.14} still holds uniformly with respect to $n \in \bbn$. We observe that in the special one-dimensional case with $A=(1,\,\infty)$, the following 
lemma can be reduced to \cite[Theorem 3.1]{tang:yuan:2016}. Further, when we say 
Assumption \ref{ass.KPT.3.1} is satisfied in following lemma, we have in mind that 
\eqref{eq.KPT.3.3} holds with common function $h(\cdot)$ for all $i \in \bbn$. 

\ble \label{lem.KPT.3.4}
Let $A \in \mathscr{R}$ be some fixed set. We consider that the $\{({\bf X}^{(i)},\,W_i)\,,\;i \in \bbn\}$ are i.i.d. copies of the pair $({\bf X},\,W)$ and satisfy Assumption \ref{ass.KPT.3.1}. If additionally we have ${\bf X} \sim F \in \mathcal{A}_A^*$ and 
\beao
\PP(0 \leq W \leq 1) =1\,,
\eeao 
$\PP(W =0) < 1$ and  $\PP(W =1) < 1$, then relation \eqref{eq.KPT.3.14} holds uniformly with respect to $n \in \bbn$, namely it holds
\beam \label{eq.KPT.3.22} 
\lim \;\sup_{n \in \bbn}\left | \dfrac{\PP\left( {\bf S}_{n} \in x\,A \right)}{ \sum_{i=1}^{n} \PP\left( {\bf X}^{(i)}\,\prod_{j=i}^{n} W_j \in x\,A\right)} - 1 \right| =0\,.
\eeam
\ele

\pr~
Let us follow the line of proof in \cite[Theorem 3.1]{tang:yuan:2016}, with the suitable modifications. At first, we show that there exists some random vector ${\bf Z}$, independent of the $\{({\bf X}^{(i)},\,W_i)\,,\;i \in \bbn\}$, such that 
\beam \label{eq.KPT.3.23} 
\PP\left( {\bf Z} \in x\,A \right) \sim c\,F( x\,A)\,,
\eeam
and the inequality 
\beam \label{eq.KPT.3.24} 
\PP\left( {\bf S}_{n} \in x\,A \right) \leq \PP\left( {\bf Z} \in x\,A\right)\,, 
\eeam
are true for any $x \geq 0$, $n \in \bbn$. By \cite[Lemma 5.3]{tang:yuan:2016}, we find 
that there exists some positive random variable $Z_A$, that is independent of $( X_A,\,W)$, with 
\beao
\PP(Z_A >x) \sim c\,\PP(X_A >x)\,,
\eeao 
for some constant $c \in (0,\,\infty)$, and for any $x \in \bbr$ it holds
\beam \label{eq.KPT.3.25} 
\PP\left( (X_A +Z_A)\,W > x \right) \leq \PP\left( Z_A > x \right)\,,
\eeam
Thus, we can find some non-negative random vector ${\bf Z}$, such that the $Z_A := \{u\;:\; {\bf Z} \in u\,A \}$ satisfies the last two conditions. Hence, relation \eqref{eq.KPT.3.23} holds 
for the random vector ${\bf Z}$. We now show relation \eqref{eq.KPT.3.24} through induction on 
$n$. For $n=1$ we have
\beao
&&\PP\left( {\bf S}_{1} \in x\,A \right) = \PP\left( X_A^{(1)} \, W_1 > x \right) \leq \PP\left( (X_A^{(1)} +Z_A)\,W_1 > x \right) \leq \PP\left( Z_A > x \right)\\[2mm]
&& =\PP\left( {\bf Z} \in x\,A\right)\,,
\eeao 
where at the third step we used relation \eqref{eq.KPT.3.25}. Let now suppose that  \eqref{eq.KPT.3.24} is true for 
some $n=k \in \bbn$. We show that it holds also for $n=k+1$. For any $x \geq 0$ we obtain
\beao
&& \PP\left( {\bf S}_{k+1} \in x\,A \right) = \PP\left( \sum_{i=1}^{k} {\bf X}^{(i)}\,\prod_{j=i}^{k+1} W_j + {\bf X}^{(k+1)}\, W_{k+1} \in x\,A \right) \\[2mm]
&& \leq \PP\left( {\bf Z}\,W_{k+1} +  {\bf X}^{(k+1)}\, W_{k+1} \in x\,A\right) \leq \PP\left( (Z_A +X_A^{(k+1)})\, W_{k+1} > x \right) \\[2mm]
&& \leq \PP\left( Z_A > x \right) = \PP\left( {\bf Z} \in x\,A\right)\,,
\eeao 
where at the second step we used the induction assumption, at the third step we took into account \cite[Proposition 2.4]{konstantinides:passalidis:2024g}, while at the 
fourth step we considered relation \eqref{eq.KPT.3.25}.

We observe that, since the $\{({\bf X}^{(i)},\,W_i)\,,\;i \in \bbn\}$ are i.i.d. copies of $({\bf X},\,W)$ we find 
\beam \label{eq.KPT.3.26} 
{\bf S}_{n} =  \sum_{i=1}^{n} {\bf X}^{(i)}\,\prod_{j=i}^{n} W_j \stackrel{d}{=}\sum_{i=1}^{n} {\bf X}^{(i)}\,\prod_{j=1}^{i} W_j =:\widehat{\bf S}_{n}\,,
\eeam    
for any $n \in \bbn$, where the '$\stackrel{d}{=}$' denotes the equality in distribution. Let us consider some fixed $M \in \bbn$. We show now that \eqref{eq.KPT.3.14} holds uniformly for $n > M$. At first we deal with the upper bound, to find that
\beam \label{eq.KPT.3.27} \notag
&& \PP\left( {\bf S}_{n} \in x\,A \right) = \PP\left( \widehat{\bf S}_{n} \in x\,A \right) \\[2mm]
&& = \PP\left( \sum_{i=1}^{M} {\bf X}^{(i)}\,\prod_{j=1}^{i} W_j +\left[ \sum_{i=M+1}^{n} {\bf X}^{(i)}\,\prod_{j=M+1}^{i} W_j \right]\,\prod_{j=1}^{M} W_j \in x\,A \right)\,,
\eeam 
for any $x \geq 0$. From relations \eqref{eq.KPT.3.26} and \eqref{eq.KPT.3.24} we obtain
\beao
&& \PP\left( \sum_{i=M+1}^{n} {\bf X}^{(i)}\,\prod_{j=M+1}^{i} W_j \in x\,A \right)= \PP\left( \sum_{i=1}^{n-M} {\bf X}^{(i)}\,\prod_{j=1}^{i} W_j \in x\,A \right) \\[2mm]
&& =\PP\left( \widehat{\bf S}_{n-M} \in x\,A \right) \leq\PP\left( {\bf Z} \in x\,A\right)\,. 
\eeao 
From \eqref{eq.KPT.3.27}, in combination with last relation we find that
\beam \label{eq.KPT.3.28} \notag
&& \PP\left( {\bf S}_{n} \in x\,A \right) = \PP\left( \widehat{\bf S}_{n} \in x\,A \right) \leq \PP\left( \sum_{i=1}^{M} {\bf X}^{(i)}\,\prod_{j=1}^{i} W_j + {\bf Z}\,\prod_{j=1}^{M} W_j \in x\,A \right)\\[2mm]
&& = \PP\left( \sum_{i=1}^{M-1} {\bf X}^{(i)}\,\prod_{j=1}^{i} W_j + \left( {\bf X}^{(M)} +{\bf Z} \right)\,\prod_{j=1}^{M} W_j \in x\,A \right)\,,
\eeam 
for any $x \geq 0$.

Next, we show that the ${\bf X}^{(M)} +{\bf Z}$, satisfy the dependence
structure in \eqref{eq.KPT.3.3} with a new function $h^*(\cdot)$, that means, it holds
\beam \label{eq.KPT.3.29} 
\PP\left( {\bf X}^{(M)} +{\bf Z} \in x\,A \;|\; W_M =y\right) \sim h^*(y)\,\PP\left({\bf X}^{(M)}+ {\bf Z} \in x\,A \right)\,,
\eeam
uniformly for any $y \in [0,\,1]$. From one hand side by \cite[Proposition 2.4]{konstantinides:passalidis:2024g} at the first step, in combination with \cite[p. 2423, after (5.7)]{tang:yuan:2016}, we have that 
\beam \label{eq.KPT.3.30} \notag
&&\PP\left( {\bf X}^{(M)} +{\bf Z} \in x\,A \;|\; W_M =y\right) \leq \PP\left( X_A^{(M)} + Z_A > x \;|\; W_M =y\right) \\[2mm]
&& \sim [h(y)+c]\,\PP\left({\bf X}^{(M)} \in x\,A \right)= [h(y)+c]\,F\left(x\,A \right)\,,
\eeam  
holds uniformly for any $y \in [0,\,1]$. 

From the other hand side, via Bonferroni inequality we obtain that it holds
\beam \label{eq.KPT.3.31} 
&&\PP\left( {\bf X}^{(M)} +{\bf Z} \in x\,A \;|\; W_M =y\right) \\[2mm] \notag
&& \geq \PP\left({\bf X}^{(M)} \in x\,A \;|\; W_M =y\right) + \PP\left( {\bf Z} \in x\,A\right)-\PP\left( {\bf Z} \in x\,A\right)\,\PP\left({\bf X}^{(M)} \in x\,A \;|\; W_M =y\right) \\[2mm] \notag
&&\sim h(y)\,\PP\left({\bf X}^{(M)} \in x\,A \right) + \PP\left( {\bf Z} \in x\,A\right) \sim  [h(y)+c]\,F\left(x\,A \right)\,,
\eeam  
uniformly for any $y \in [0,\,1]$, where at the last step we used \eqref{eq.KPT.3.23}. From the fact that $X_A^{M},\,Z_A$ have distributions from the class $\mathcal{S}$, 
\beao
\PP\left(  Z_A > x\right) \sim c\,\PP\left( X_A^{(M)} > x\right)\,,
\eeao 
with $c \in (0,\,\infty)$, and $X_A^{(M)},\,Z_A$ independent each other, is implied that  the sum $X_A^{(M)} + Z_A$ 
follows distribution from the class $\mathcal{S}$, due to \cite[Lemma 3.2]{tang:tsitsiashvili:2003a}. This, in combination with 
\cite[Theorem 3.5]{konstantinides:passalidis:2024g} is implied that the sum ${\bf X}^{(M)} +{\bf Z}$, follows distribution from the class $\mathcal{S}_A$, and further we find
\beam \label{eq.KPT.3.32} 
\PP\left( {\bf X}^{(M)} +{\bf Z} \in x\,A \right) \sim \PP\left({\bf X}^{(M)} \in x\,A \right) + \PP\left( {\bf Z} \in x\,A\right) \sim (c+1)\,F\left(x\,A \right)\,.
\eeam
Hence, from \eqref{eq.KPT.3.30}, \eqref{eq.KPT.3.31}, in combination with \eqref{eq.KPT.3.32}, we obtain that \eqref{eq.KPT.3.29} is true with 
\beao
h^*(y) = \dfrac{h(y)+c}{c+1}\,.
\eeao                                          
Therefore, from \eqref{eq.KPT.3.28} we find that
\beam \label{eq.KPT.3.33} \notag
&& \PP\left( {\bf S}_{n} \in x\,A \right) = \PP\left( \widehat{\bf S}_{n} \in x\,A \right)  \\[2mm] \notag
&&\leq \PP\left( \sum_{i=1}^{M-1} {\bf X}^{(i)}\,\prod_{j=1}^{i} W_j +\left( {\bf X}^{(M)} +{\bf Z} \right)\,\prod_{j=1}^{M} W_j \in x\,A \right)\\[2mm]
&& \sim  \sum_{i=1}^{M-1} \PP\left( {\bf X}^{(i)}\,\prod_{j=1}^{i} W_j \in x\,A \right) + \PP\left(\,({\bf X}^{(M)} +  {\bf Z} )\,\prod_{j=1}^{M} W_j \in x\,A \right) \\[2mm] \notag
&& \sim  \sum_{i=1}^{M-1} \PP\left( {\bf X}^{(i)}\,\prod_{j=1}^{i} W_j \in x\,A \right) + \PP\left({\bf X}^{(M)}\,\prod_{j=1}^{M} W_j \in x\,A \right) + \PP\left( {\bf Z}\,\prod_{j=1}^{M} W_j \in x\,A \right) \\[2mm] \notag
&& \leq \sum_{i=1}^{n} \PP\left( {\bf X}^{(i)}\,\prod_{j=1}^{i} W_j \in x\,A \right) + \PP\left( {\bf Z}\,\prod_{j=1}^{M} W_j \in x\,A \right)\,, 
\eeam
where at the third step we used Lemma \ref{lem.KPT.3.3}; recall relation \eqref{eq.KPT.3.26} in the first term, and relations \eqref{eq.KPT.3.26} and  \eqref{eq.KPT.3.32} in the second term, while at the fourth step we applied Lemma \ref{lem.KPT.3.1}.

For the second term in last relation of \eqref{eq.KPT.3.33}, by \cite[Lemma 5.2]{tang:yuan:2016}, we obtain that for any 
\beao
\delta \in (0,\,K_{F_A}^-)\,,
\eeao 
recall that, since 
\beao
\PP\left( Z_A > x\right) \sim c\,\bF_A\left( x\right)\,,
\eeao 
the $ Z_A$ distribution has the same lower Karamata index with that of $F_A$, and for any
\beao
a>1\,,
\eeao 
there exists some $x_0> 0$ independent of $M$, such that it holds
\beam \label{eq.KPT.3.34} \notag
&& \PP\left( {\bf Z}\,\prod_{j=1}^{M} W_j \in x\,A \right) = \PP\left( Z_{A}\,\prod_{j=1}^{M} W_j >x \right) =\int_0^1 \PP\left( Z_{A}\,\prod_{j=2}^{M} W_j >\dfrac xy \right)\,\PP\left( W_1 \in dy\right)  \\[2mm] \notag
&&\leq a\,\E\left[ \prod_{j=2}^{M} W_j^{K_{F_A}^- - \delta} \right]\,\int_0^1 \PP\left( Z_{A} >\dfrac xy \right)\,\PP\left( W_1 \in dy\right) =  a\,\left(\E\left[  W^{K_{F_A}^- - \delta} \right] \right)^{M-1}\\[2mm]
&&\times \,\PP\left( {\bf Z}\, W \in x\,A \right) \leq C\,a\,\left(\E\left[  W^{K_{F_A}^- - \delta} \right] \right)^{M-1}\,\PP\left( {\bf X}\,W \in x\,A \right)\,, 
\eeam 
for any $x >x_0$, where the constant $C>0$, follows from the fact that 
\beao
\PP\left( {\bf Z}\,W \in x\,A \right) \asymp \PP\left( {\bf X}\,W \in x\,A \right)\,.
\eeao 
Indeed, from the fact that the function $h(\cdot)$ is bounded away from zero and infinity on a left neighborhood of unity, say $E$, we obtain that
\beao
&&\PP\left( {\bf X}\,W \in x\,A \right) = \int_0^1 \PP\left( {\bf X}\in \dfrac xy \,A \;\big|\; W =y\right)\,\PP\left( W \in dy\right)  \\[2mm]
&& \sim \int_0^1 h(y )\,\PP\left( {\bf X}\in \dfrac xy \,A \right)\,\PP\left( W \in dy\right)  \asymp  \int_{E} h(y)\, \PP\left( {\bf X}\in \dfrac xy \,A\right)\,\PP\left( W \in dy\right)\\[2mm]
&& \asymp \int_{E} \PP\left( {\bf X}\in \dfrac xy \,A\right)\,\PP\left( W \in dy\right) \sim  c\,\int_{E} \PP\left( {\bf Z}\in \dfrac xy \,A\right)\,\PP\left( W \in dy\right)\\[2mm]
&& \asymp \int_0^1 \PP\left( {\bf Z}\in \dfrac xy \,A\right)\,\PP\left( W \in dy\right) =  \PP\left( {\bf Z}\,W \in x\,A \right) \,.
\eeao 
Hence, since 
\beao
\E\left[ Y^{K_{F_A}^- - \delta} \right] <1\,,
\eeao 
for any fixed $\vep >0$, we can find a sufficiently large $M \in \bbn$, such that we obtain
~\halmos
\beao
\PP\left(  {\bf Z}\,\prod_{j=1}^{M} W_j \in x\,A\right) \leq \vep\,\PP\left( {\bf X}\,W \in x\,A \right)\,,
\eeao
for any $x > x_0$; recall relation \eqref{eq.KPT.3.34}. Thus, from this last inequality, in combination with \eqref{eq.KPT.3.33} we have that
\beam \label{eq.KPT.3.35} 
&& \PP\left( {\bf S}_{n} \in x\,A \right) = \PP\left( \widehat{\bf S}_{n} \in x\,A \right)  \\[2mm] \notag
&& \lesssim (1+ \vep)\,\sum_{i=1}^{n} \PP\left( {\bf X}^{(i)}\,\prod_{j=1}^{i} W_j \in x\,A \right)= (1+ \vep)\,\sum_{i=1}^{n} \PP\left( {\bf X}^{(i)}\,\prod_{j=i}^{n} W_j \in x\,A \right)\,,
\eeam 
where at the last equality we used the fact that the $\{({\bf X}^{(i)},\,W_i)\,,\;i \in \bbn\}$ are i.i.d. copies of the pair $({\bf X},\,W)$.

Now we deal with the lower bound of \eqref{eq.KPT.3.14} uniformly for $n > M$. We start via Lemma \ref{lem.KPT.3.3}, by obtaining the following
\beam \label{eq.KPT.3.36}
&& \PP\left( {\bf S}_{n} \in x\,A \right) \geq \PP\left( {\bf S}_{M} \in x\,A \right)\\[2mm] \notag
&& \sim \sum_{i=1}^{M} \PP\left( {\bf X}^{(i)}\,\prod_{j=i}^{M} W_j \in x\,A \right)  = \sum_{i=1}^{M} \PP\left( {\bf X}^{(i)}\,\prod_{j=1}^{i} W_j \in x\,A \right)\\[2mm] \notag
&&= \left(\sum_{i=1}^{n} - \sum_{i=M+1}^{n} \right)\,\PP\left( {\bf X}^{(i)}\,\prod_{j=1}^{i} W_j \in x\,A \right)\,.
\eeam
Hence, again due to \cite[Lemma 5.2]{tang:yuan:2016} we find for any $\delta \in (0,\,K_{F_A}^-)$ 
and any $a>1$, that there exists some $x_0>0$, independent of $M$, such that for any $x > x_0$ it holds
\beao
&&  \sum_{i=M+1}^{n} \PP\left( {\bf X}^{(i)}\,\prod_{j=1}^{i} W_j \in x\,A \right)\\[2mm] 
&&=\sum_{i=M+1}^{n} \int_0^1 \int_0^1 \PP\left( {\bf X}^{(i)} \in \dfrac x{y\,s} \,A \;\big|\; W_i =y\right)\,\PP\left( \prod_{j=1}^{i-1} W_j \in ds \right)\,\PP\left( W_i \in dy\right)\\[2mm] 
&&\sim \sum_{i=M+1}^{n} \int_0^1  h(y)\,\PP\left( {\bf X}^{(i)}\,\prod_{j=1}^{i-1} W_j \in \dfrac x{y} \,A \right)\,\PP\left( W_i \in dy\right)\\[2mm] 
&&\leq \sum_{i=M+1}^{n} a\,\E\left[ \prod_{j=1}^{i-1} W_j^{K_{F_A}^- - \delta} \right]\, \int_0^1 h(y)\,\PP\left( {\bf X}^{(i)} \in \dfrac x{y} \,A \right)\,\PP\left( W_i \in dy\right)\\[2mm] 
&&\sim \sum_{i=M+1}^{n} a\,\left(\E\left[ W^{K_{F_A}^- - \delta} \right]\right)^{i-1}\,\PP\left( {\bf X} \,W\in x\,A \right) \leq a\,\dfrac{\left(\E\left[ W^{K_{F_A}^- - \delta} \right] \right)^{M}}{1-\E\left[ W^{K_{F_A}^- - \delta} \right]}\,\PP\left( {\bf X} \,W\in x\,A \right)\,,
\eeao 
where at the next-to-last step we used \eqref{eq.KPT.3.8}. Hence, since $\E\left[ W^{K_{F_A}^- - \delta} \right] <1$, from the last relation we can find, for any $\vep >0$, some sufficiently large $M \in \bbn$, such that 
\beao
\sum_{i=M+1}^{n} \PP\left( {\bf X}^{(i)}\,\prod_{j=1}^{i} W_j \in x\,A \right) \lesssim \vep\,\PP\left( {\bf X} \,W\in x\,A \right)\,,
\eeao
and further from this last relation, in combination with relation \eqref{eq.KPT.3.36} we obtain that for any $n> M$
\beam \label{eq.KPT.3.37}
&& \PP\left( {\bf S}_{n} \in x\,A \right) \gtrsim (1-\vep)\,\sum_{i=1}^{n} \PP\left( {\bf X}^{(i)}\,\prod_{j=1}^{i} W_j \in x\,A \right) \\[2mm] \notag
&& = (1-\vep)\,\sum_{i=1}^{n} \PP\left( {\bf X}^{(i)}\,\prod_{j=i}^{n} W_j \in x\,A \right)\,.
\eeam
So, from relations \eqref{eq.KPT.3.35} and  \eqref{eq.KPT.3.37}, we have that relation  \eqref{eq.KPT.3.14} holds uniformly for any $n>M$. Relation \eqref{eq.KPT.3.14} holds uniformly for $1 \leq n \leq M$, which follows by Lemma \ref{lem.KPT.3.3}.
~\halmos

\bre \label{rem.KPT.a}
Let us notice that under the conditions of Lemma \ref{lem.KPT.3.4}, we can see that we have
\beam \label{eq.KPT.3.38}
\PP\left( \widehat{\bf S}_{n} \in x\,A \right) \sim \sum_{i=1}^{n} \PP\left( {\bf X}^{(i)}\,\prod_{j=1}^{i} W_j \in x\,A \right)\,,
\eeam
uniformly for $n \in \bbn$; recall the relations \eqref{eq.KPT.3.35} and \eqref{eq.KPT.3.37}, in combination with \eqref{eq.KPT.3.26}.
\ere

Now we are ready to present the proof of Theorem \ref{th.KPT.3.1}

{\bf Proof of Theorem \ref{th.KPT.3.1}}~
Since the $\{R(t)\,,\;t\geq 0\}$ is non-negative and non-degenerative to zero L\'{e}vy process, we find that the relations 
\beam \label{eq.KPT.3.39}
W_i := e^{R(\tau_{i-1}) - R(\tau_{i})}\stackrel{d}{=} e^{-R(\theta_{1})}\,,
\eeam
for any $i \in \bbn$, are i.i.d., and it holds 
\beao
\PP(0\leq W_i \leq 1) =1\,,
\eeao 
$\PP(W_i =0)< 1$ and $\PP(W_i =1) < 1$. We also have that
\beam \label{eq.KPT.3.40}
\PP\left({\bf X} \in x\,A \;|\; e^{-R(\theta_{1})} =y\right) \sim h(y)\,\PP\left({\bf X} \in x\,A \right)\,,
\eeam
uniformly for any $y \in [0,\,1]$, due to Assumption \ref{ass.KPT.2.2}. Indeed, we obtain
\beao
&&\PP\left({\bf X} \in x\,A \;|\; e^{R(-\theta_{1})} =y\right) = \int_{\Lambda} \dfrac{\PP\left({\bf X} \in x\,A \,,\; e^{-R(s_{1})} =y\;|\; \theta_{1} =s_1\right)}{\PP\left(e^{-R(s_{1})} =y\;|\; \theta_{1} =s_1\right)}\,\PP(\theta_1 \in ds_1) \\[2mm]
&&= \int_{\Lambda} \dfrac{\PP\left({\bf X} \in x\,A \,,\; e^{-R(s_{1})} =y\,,\; \theta_{1} =s_1\right)/\PP(\theta_1 \in s_1)}{\PP\left(e^{-R(s_{1})} =y\,,\; \theta_{1} =s_1\right)/\PP(\theta_1 \in s_1)\,}\,\PP(\theta_1 \in ds_1) \\[2mm]
&& = \int_{\Lambda} \PP\left({\bf X} \in x\,A \;|\; e^{-R(s_{1})} =y\,,\; \theta_{1} =s_1\right)\,\PP(\theta_1 \in ds_1) \sim \int_{\Lambda} h(y)\,\PP\left({\bf X} \in x\,A \right)\,\PP(\theta_1 \in ds_1)\\[2mm]
&&= h(y)\,\PP\left({\bf X} \in x\,A \right)\,,
\eeao
where at the next-to-last step we used the dominated convergence theorem, because of the uniformity of \eqref{eq.KPT.2.8} with respect to $s_1 \in \Lambda$. The uniformity of \eqref{eq.KPT.3.40} with respect to $y \in [0,\,1]$, remains also true at the next-to-last step in last relation, via relation \eqref{eq.KPT.2.8}. Hence, by \eqref{eq.KPT.3.39} and \eqref{eq.KPT.3.40}, applying Lemma \ref{lem.KPT.3.4}, in the form of \eqref{eq.KPT.3.38}, for $n=\infty$, we finally find
\beao
&&\PP\left({\bf D}(\infty) \in x\,A \right) = \PP\left(\sum_{i=1}^{\infty} {\bf X}^{(i)}\,e^{-R(\tau_{i})} \in x\,A\right) \\[2mm]
&&\sim \sum_{i=1}^{\infty} \PP\left( {\bf X}^{(i)}\,e^{-R(\tau_{i})} \in x\,A\right) = \int_0^{\infty} \PP\left({\bf X}\,e^{-R(s)} \in x\,A \right)\,\lambda (ds)\,.~\halmos
\eeao

\section{Arbitrarily dependent insurance and financial Risks}

\subsection{Main result } \label{sec.KLP.4.1}

In this section  we look for asymptotic behavior of the entrance probability of discounted aggregate claim vectors into some rare set $x\,A$, on infinite time horizon, when insurance and financial risks are arbitrarily dependent. In opposite to previous section, now we assume  that the condition $H \in (\mathcal{D}\cap \mathcal{A})_A$ holds, where
	\beam \label{eq.KLP.3.1}
	H(x\,A) = \PP( {\bf X}\,e^{-R(\theta_1)} \in x\,A) =\PP(X_A\,e^{-R(\theta_1)} > x ) = \bH_A(x)\,,
	\eeam 
and now we have remove the condition of non-negative L\'{e}vy process and instead of this we assume that the Laplace exponent satisfies $\phi(p)< 0$, for some $p>J_{H_A}^+$. Namely, the distribution $H$ depicts the product distribution of the insurance and financial risks, at each renewal epoch. 
	
The following assumption is crucial in all over the section.
\begin{assumption} \label{ass.KLP27.1.1}
We suppose that $\{{\bf X}^{(i)}\,e^{R(\tau_{i-1})-R(\tau_i)}\,,\;i \in \bbn\}$ represents a sequence of i.i.d. random vectors with common distribution $H$.
\end{assumption}

Now we are ready to provide the main theorem, in which the insurance and financial risks are arbitrarily dependent.

\bth \label{th.KLP.4.1}
Let  $A \in \mathscr{R}$ be a fixed set and let consider the discounted aggregate claim vectors in relation \eqref{eq.KLP.1.3}. Under the Assumptions \ref{ass.KPT.2.1} and \ref{ass.KLP27.1.1}, with $H \in (\mathcal{D}\cap \mathcal{A})_A $ and additionally if $\phi(p)<0$, for some $p>J_{H_A}^+$, then we obtain that
\beam \label{eq.KLP.4.1} 
\PP\left[{\bf D}(\infty) \in x\,A\right]\sim \int_0^{\infty} \PP\left[{\bf X}\,e^{-R(s)} \in x\,A\right]\,\lambda(ds)\,.
\eeam
\ethe

\bre \label{rem.KLP.4.1}
As in Theorem \ref{th.KPT.3.1}, if here we demand that $H \in (\mathcal{D}\cap \mathcal{A})_\mathscr{R} $ and $\phi(p)<0$, for some $p>J_{H_A}^+$, for any $A \in \mathscr{R}$, then relation \eqref{eq.KLP.4.1} holds for any $A \in \mathscr{R}$. Even more, in special case when $A=(1,\,\infty)$ the result of Theorem \ref{th.KLP.4.1} is new. In comparison with \cite[Th. 2]{cheng:konstantinides:wang:2024}, in which was considered $d-$L\'{e}vy processes, that were arbitrarily dependent, hence it was more general model, was estimated the ruin probability over infinite time horizon, under the condition of $MRV$ insurance and financial risks, with asymptotically dependent components, namely $\mu(({\bf 1},\,\vec{\infty}])>0$, and similar assumptions on Laplace exponent of the L\'{e}vy processes. Here, although we use a common L\'{e}vy process, the asymptotic estimation for the distribution tail of the discounted aggregate claim vectors is under the distribution class $ (\mathcal{D}\cap \mathcal{A})_A$, which is larger than $MRV$, and we do NOT demand asymptotic dependence among the components, see in \cite[Exam. 4.2, 4.3]{konstantinides:liu:passalidis:2025}, for two  distributions that belong in $ (\mathcal{D}\cap \mathcal{A})_A$, with asymptotically independent components.
\ere

\bre \label{rem.KLP.4.2}
As we shall see in the subsection \ref{sec.KLP.5}, the framework of Theorem \ref{th.KLP.4.1} allows a wide type of dependencies between insurance and financial risks, and contains both asymptotic independent and asymptotic dependent cases. This framework allows also the domination of financial risks into the insurance risks, instead of the classical problems in risk theory in which only the insurance risks dominated in financial risks.
\ere

\subsection{Discussion on dependence between the risks}\label{sec.KLP.5}

In this section we provide two examples, that satisfy the conditions a)$H \in (\mathcal{D}\cap \mathcal{A})_A$ and b) $\phi(p)< 0$, for some $p > J_{H_A}^+$. These examples are based to conditions  only on the marginal distributions of the insurance and financial risks, and also the dependence structures between these two risks. This way, the assumptions of the Theorem \ref{th.KLP.4.1} are easily controllable, but also represent a usable tool in practice. In order to stress on the importance of the dependence between insurance and financial risks, we give an example with the weak dependence structure of Assumption \ref{ass.KPT.2.2}, and another one with a strong dependence structure. 
	
\bexam \label{exam.KLP.5.1}
Let $A \in \mathscr{R}$ be some fixed set, and let Assumptions \ref{ass.KPT.2.1} and \ref{ass.KPT.2.2} holds. If $F \in (\mathcal{D}\cap \mathcal{A})_A$ and $\{R(t)\,,\;t\geq 0\}$ is a L\'{e}vy process with $\phi(p)< 0$, for some $p > J_{F_A}^+$, then we obtain $H \in (\mathcal{D}\cap \mathcal{A})_A$ and $\phi(p^*)< 0$, for some $p* > J_{H_A}^+$.
\eexam
	
	\pr~
We show that $H \in (\mathcal{D}\cap \mathcal{A})_A$, that follows by
$H_A \in \mathcal{D}\cap \mathcal{A}$, recall relation \eqref{eq.KLP.3.1}. Since 
	$\phi(p)< 0$ for some $p> J_{F_A}^+$, we find that
\beam \label{eq.KLP.5.3}
\E\left[ e^{-p\,R(\theta_{1})} \right] =\E\left[ e^{\theta_{1}\,\phi(p)} \right] < 1
\eeam
through \eqref{eq.KLP.1.1}. Hence we obtain that
\beam \label{eq.KLP.5.4}
&&\bH_A(x)=\PP\left[ X_A\,e^{-R(\theta_{1})}>x \right] =\int_{\Lambda}\PP\left[ 
X_A\,e^{-R(s)}>x \right]\,\PP(\theta_1 \in ds)\\[2mm]\notag
&& =\int_{\Lambda} \int_0^{\infty} \PP\left[X_A> \dfrac xy\;|\;e^{-R(s)} =y \right]
\PP\left[ e^{-R(s)} \in dy \right] \PP(\theta_1 \in ds)=: \int_{\Lambda} I(x,\,s)\PP(\theta_1 \in ds).
\eeam
Let consider some constant $k$ such that it satisfies $J_{F_A}^+/p < k <1$
and we separate $I(x,\,s)$ into two parts as follows
\beao
&&I(x,\,s)=\left(\int_0^{x^k} + \int_{x^k}^{\infty} \right) \PP\left[X_A> \dfrac xy\;|\;e^{-R(s)} =y \right]\,\PP\left[ e^{-R(s)} \in dy \right]\\[2mm]
&&= I_1(x,\,s,\,k)+ I_2(x,\,s,\,k)\,.
\eeao 
Integrating with respect to the distribution of $\theta_1$, the second term, and applying Markov inequality we obtain that
\beam \label{eq.KLP.5.6} \notag
&&\int_{\Lambda} I_2(x,\,s, \,k)\,\PP(\theta_1 \in ds) \leq \int_{\Lambda} \PP\left[e^{-R(s)} > x^k \right]\,\PP(\theta_1 \in ds) =\\[2mm]
&&=\PP\left[e^{-R(\theta_1)} >x^k \right] \leq \E\left[ e^{-p\,R(\theta_1)} \right]\,x^{-p\,k}=o(\PP(X_A>x)\,)\,,
\eeam
where we used relation \eqref{eq.KLP.5.3} and the fact that $p\,k> J_{F_A}^+$, which allows us to apply \eqref{eq.KP.2.3a}. 
		
Further we deal with $I_1(x,\,s,\,k)$. Let define a new c\'{a}dl\'{a}g process $\{R^*(t)\,,\;t\geq 0\}$, which is independent of all the other sources of randomness, which satisfies the relation \eqref{eq.KPT.2.10}. Considering the c\'{a}dl\'{a}g process $\{R^*(t)\,,\;t\geq 0\}$, we find that on $p$ it holds
\beao
&&\E\left[ e^{-p\,R^*(\theta_1)} \right]= \int_{\Lambda} \E\left[ e^{-p\,R^*(s)} \right]\,\PP(\theta_1 \in ds) =\int_{\Lambda} \int_0^{\infty} y^p\,\PP\left[ e^{-R^*(s)} \in dy\right]\,\PP(\theta_1 \in ds) \\[2mm]
&&=\int_{\Lambda} \int_0^{\infty} y^p\,h(y)\,\PP\left[ e^{-R(s)} \in dy\right]\,\PP(\theta_1 \in ds) =\int_{\Lambda} \E\left[e^{-p\,R(s)}\,h\left( e^{-R(s)}\right)\right]\,\PP(\theta_1 \in ds) \\[2mm]
&&=\E\left[e^{-p\,R(\theta_{1})}\,h\left( e^{-R(\theta_{1})}\right) \right] \leq c\,\E\left[e^{-p\,R(\theta_1)} \right]\leq c\,<\infty\,,
\eeao  
where at the third step we used relation \eqref{eq.KPT.2.10}, at the penultimate step the constant $c>0$ comes from the fact that function $h(y)$ is bounded from above, see in \cite[Prop. 2.4]{cui:wang:2025} when $y$ represents a probable value, while it takes the value $c=1$ in case when $y$ is not a probable value. Hence, we obtain that it holds
\beam \label{eq.KLP.5.8} 
\E\left[e^{-p\,R^*(\theta_1)} \right] <\infty\,,
\eeam
and since the $e^{-R^*(\theta_1)} $ is independent of the random variable $X_A$, whose distribution is $F_A \in \mathcal{D}\cap \mathcal{A}$, in view of \cite[Th. 3.3(iv)]{cline:samorodnitsky:1994} and \cite[Th. 5.1(i)]{konstantinides:passalidis:2025b}, we obtain that the product $X_A\,e^{-R^*(\theta_1)}$ follows distribution $H_A^*\in \mathcal{D}\cap \mathcal{A}$, and further again by \cite[Th. 3.3(iv)]{cline:samorodnitsky:1994} it holds
\beam \label{eq.KLP.5.9} 
\PP\left[X_A\,e^{-R^*(\theta_1)} >x \right] \asymp \PP\left[X_A >x \right] \,.
\eeam
Therefore, for the estimation of $I_1(x,\,s,\,k)$ we obtain that
\beam \label{eq.KLP.5.10} \notag
&&\int_{\Lambda} I_1(x,\,s,\,k)\,\PP(\theta_1 \in ds) =\int_{\Lambda} \int_0^{x^k} \PP\left[X_A>\dfrac xy\;|\; e^{-R(s)} = y\right]\,\PP\left[ e^{-R(s)} \in dy\right]\,\PP(\theta_1 \in ds) \\[2mm] \notag
&& \sim\int_{\Lambda} \int_0^{x^k}  \,h(y)\PP\left[X_A>\dfrac xy\right]\,\PP\left[ e^{-R(s)} \in dy\right] \,\PP(\theta_1 \in ds) \\[2mm] \notag
&&=\int_{\Lambda} \left( \int_0^{\infty} - \int_{x^k}^{\infty}\right) \,\PP\left[X_A>\dfrac xy\right]\,\PP\left[ e^{-R^*(s)} \in dy\right] \,\PP(\theta_1 \in ds) \\[2mm]
&&=\PP\left[X_A\, e^{-R^*(\theta_1)} >x\right]- o(\PP(X_A>x)\,)  \sim \PP\left[X_A\, e^{-R^*(\theta_1)} >x\right]\,,
\eeam
where at the second step we used relation \eqref{eq.KPT.2.8}, due to dominated convergence theorem, recall that $k<1$ and the uniformity of \eqref{eq.KPT.2.8}. Further at the third step we used relation \eqref{eq.KPT.2.10}, while at the fourth step the second term is negligible, in view of the relation \eqref{eq.KLP.5.6}, keeping in mind that the function $h(\cdot)$ is bounded from above. At the last step we used relation  \eqref{eq.KLP.5.9}. Hence, from relations  \eqref{eq.KLP.5.6}, \eqref{eq.KLP.5.9} and \eqref{eq.KLP.5.10}, in combination with relation \eqref{eq.KLP.5.4} we obtain that it holds
\beam \label{eq.KLP.5.11} 
\bH_A(x)=\PP\left[X_A\, e^{-R(\theta_1)} >x\right] \sim \PP\left[X_A\, e^{-R^*(\theta_1)} >x\right]=\bH_A^*(x)\,.
\eeam
Since $H_A^* \in \mathcal{D}\cap \mathcal{A}$, by relation \eqref{eq.KLP.5.11} we find that $H_A \in \mathcal{D}\cap \mathcal{A}$. Furthermore, from \cite[Lem. 3.9)]{tang:tsitsiashvili:2003a} and relation \eqref{eq.KLP.5.8}, we obtain that it holds $J_{H_A^*}^+ = J_{F_A}^+$, while from relation \eqref{eq.KLP.5.11} we find that $J_{H_A^*}^+ = J_{H_A}^+$. Thus, we find $\phi(p)<0$ for some $p > J_{H_A}^+$. 
~\halmos

The following Corollary gives a more explicit form of relation \eqref{eq.KLP.4.1}, under the assumptions of Example \ref{exam.KLP.5.1}, when $F\in MRV$. Recall that $G \in \mathcal{R}_{-\alpha}$, with $a \in (0,\,\infty)$, is the auxiliary distribution function in $MRV$.
	
\bco \label{cor.KLP.5.1}
Let $A \in \mathscr{R}$ be some fixed set. Then	under the conditions of Example \ref{exam.KLP.5.1}, under the restriction 
$F \in MRV(\alpha,\,\mu)$, with $\alpha \in (0,\,\infty)$, it holds
\beam \label{eq.KLP.5.13} 
&&\PP\left[{\bf D}(\infty) \in x\,A\right] \sim\\[2mm] \notag
&& \mu(A) \bG(x) \E\left[e^{-\alpha\,R(\theta_1)} h\left(e^{-R(\theta_1)} \right)\right] \left(1 +\int_{\Lambda} \int_0^{\infty}  \E\left[e^{-\alpha R(s)} \right]\,\lambda(ds)\,\PP(\theta_1 \in du) \right).
\eeam
\eco
	
\pr~
We initially observe that the combination of Assumptions \ref{ass.KPT.2.1} and \ref{ass.KPT.2.2}, with the fact that $\{R(t)\,,\; t\geq 0\}$ is L\'{e}vy process, imply that Assumption \ref{ass.KLP27.1.1} holds. Hence, from Theorem \ref{th.KLP.4.1}, see also Example \ref{exam.KLP.5.1}, it is enough to show that
\beam \label{eq.KLP.alpha}
&&\int_0^{\infty} \PP\left[ {\bf X}\,e^{-R(s)} \in x\,A \right]\,\lambda(ds) \sim\\[2mm] \notag
&& \mu(A)\,\bG(x)\,\E\left[e^{-\alpha\,R(\theta_1)}\,h\left(e^{-R(\theta_1)} \right)\right] \left(1 +\int_{\Lambda} \int_0^{\infty}  \E\left[e^{-\alpha\,R(s)} \right]\,\lambda(ds)\,\PP(\theta_1 \in du)\, \right)\,.
\eeam
At first, since $F \in MRV(\alpha,\,\mu)$, for any $A \in \mathscr{R}$ it is true that $F_A \in \mathcal{R}_{-\alpha}$, while 
due to relation \eqref{eq.KLP.5.3}, we can apply Breiman's theorem at the third step and at forth set due to \cite[Lem. 3.3]{konstantinides:leipus:passalidis:siaulys:2024} we obtain that $X_A\,e^{-R(\theta_1)} \sim H_A\in \mathcal{R}_{-\alpha}$. Hence, 
\beao
&&\int_0^{\infty} \PP\left[ {\bf X}\,e^{-R(s)} \in x\,A \right]\,\lambda(ds) = \sum_{i=1}^{\infty} \PP\left[{\bf X}^{(i)}\,e^{-R(\tau_i)} \in x\,A \right] \\[2mm]
&& =\sum_{i=1}^{\infty} \PP\left[ X_{A}^{(i)}\,e^{R(\tau_{i-1})-R(\tau_i)}\,e^{-R(\tau_{i-1})} \in x\,A \right] \\[2mm]
&& \sim \sum_{i=1}^{\infty} \E\left[e^{-\alpha\,R(\tau_{i-1})} \right]\,\PP\left[ X_{A}^{(i)}\,e^{R(\tau_{i-1})-R(\tau_i)}>x \right] \\[2mm]
&& \sim \PP\left[ {\bf X} \in x\,A \right]\,\sum_{i=1}^{\infty} \E\left[e^{-\alpha\,R(\tau_{i-1})} \right]\,\E\left[ e^{\alpha[R(\tau_{i-1})-R(\tau_i)]}\,h\left(e^{R(\tau_{i-1})-R(\tau_i)} \right) \right] \\[2mm]
&& \sim \mu(A)\,\bG(x)\,\E\left[e^{-\alpha\,R(\theta_1)}\,h\left(e^{-R(\theta_1)} \right)\right]\,\sum_{i=1}^{\infty} \E\left[e^{-\alpha\,R(\tau_{i-1})} \right]\, \\[2mm]
&& =\mu(A)\,\bG(x)\,\E\left[e^{-\alpha\,R(\theta_1)}\,h\left(e^{-R(\theta_1)} \right)\right]\,\sum_{j=0}^{\infty} \E\left[e^{-\alpha\,R(\tau_{j})} \,{\bf 1}_{\{\tau_{j+1}< \infty\}}\right] = \mu(A)\,\\[2mm]
&&\times \bG(x)\,\E\left[e^{-\alpha\,R(\theta_1)}\,h\left(e^{-R(\theta_1)} \right)\right]\,\left(1 + \sum_{j=1}^{\infty} \int_{\Lambda}  \E\left[e^{-\alpha\,R(\tau_j)} \right]\,\PP(\tau_j < \infty - u)\,\PP(\theta_1 \in du)\, \right)\\[2mm]
&& \sim \mu(A)\,\bG(x)\,\E\left[e^{-\alpha\,R(\theta_1)}\,h\left(e^{-R(\theta_1)} \right)\right]\,\left(1 +\int_{\Lambda} \int_0^{\infty}  \E\left[e^{-\alpha\,R(s)} \right]\,\lambda(ds)\,\PP(\theta_1 \in du)\, \right)\,,	
\eeao
where at the sixth step we made change of summation variables, while we mention that $\lambda(t) < \infty$, for $t>0$, hence the $\{\theta_i\,,\;i \in \bbn\}$ have NOT mass at infinity.
~\halmos
\bre \label{rem.KLP.5.1}
If follows easily, that the dependence between insurance and financial risks, is represented entirely only through the function $h(\cdot)$ in relation \eqref{eq.KLP.5.13}, in the presence of $MRV$. This means that the relation \eqref{eq.KLP.5.13} provide an easy enough tool in practical applications. Unfortunately, in the cases where $F \notin MRV$ the corresponding relations \eqref{eq.KLP.1.4} and \eqref{eq.KLP.4.1} are much more difficult with respect to calculations. We also note that the asymptotic relation \eqref{eq.KLP.5.13} is sharp, since for any $A \in \mathscr{R}$, we obtain $\mu(A) \in (0,\,\infty)$, see for this in the proof of \cite[Proposition 4.14]{samorodnitsky:sun:2016}. 
\ere
	
\bre \label{rem.KLP.5.2}
In Example \ref{exam.KLP.5.1} the insurance risks dominate of the financial risks, in the sense that they have 	heavier tails. This can be checked by relation \eqref{eq.KLP.5.3}, due to condition $\phi(p) <0$, for some $p > J_{F_A}^+$. Intuitively we believe that, in any case of weak dependence, that means under asymptotic independence, the results of section \ref{sec.KLP.4.1} are restricted in the case where the insurance risks dominate on the financial risks, because of the Laplace exponent conditions. Although most of the papers on risk theory, consider domination of the insurance risks on the financial ones, the opposite case is interesting 	too. In discrete time risk models, we find some papers on this direction, as \cite{li:tang:2015}, under independence between these two kind of risks, and we find a generalization with a weak dependence structure (Sarmanov), in \cite{yang:jiang:wang:yuen:2020}.
\ere
	
Although, in Example \ref{exam.KLP.5.1} the insurance risks dominate on the financial risks and there exists weak dependence between the two kinds of risk, in next example there exists strong dependence structure, that means asymptotic dependence, without some domination assumption among risks. The interesting feature of the following example is that the product of the insurance and financial risks has heavier tails than any of the two risks. 
	
\bexam \label{exam.KLP.5.2}
Let $A \in \mathscr{R}$ be some fixed set. Let assume that $e^{-R(\theta_1)}$ follows distribution $Q$, which belongs in class $\mathcal{R}_{-\beta}$, with $\beta \in (0,\,\infty)$, and that $F \in MRV(\alpha,\,\mu)$ with $\alpha \in (0,\,\infty)$, and we write it, in equivalent form with normalization function $b(\cdot)$ as follows
\beao
\lim x\,\PP\left[\dfrac{{\bf X}}{b(x)} \in \bbb \right] = \mu(\bbb)\,,
\eeao
for any $\mu$-continuous Borel set $\bbb \subsetneq [0,\,\infty]^d \setminus \{{\bf 0}\}$. If we denote by $c(\cdot)$ the normalization function of $e^{-R(\theta_1)}$, then we assume that the pair $(F,\,Q)$ satisfies the non-standard $MRV$, namely we have 
\beam \label{eq.KLP.5.15} 
\lim x\,\PP\left[\left(\dfrac{{\bf X}}{b(x)}\,,\; \dfrac{e^{-R(\theta_1)}}{c(x)}\right)\in \bbb \right]=\mu^*(\bbb)\,,
\eeam
(with $\mu^*$ a new, non-degenerate to zero, Radon measure) for any $\mu^*$-continuous Borel set $\bbb \subsetneq [0,\,\infty]^{d+1} \setminus \{{\bf 0}\}$. 	Furthermore, we assume that for some $j=1,\,\ldots,\,d$ it holds
\beam \label{eq.KLP.5.16} 
\lim x\,\PP\left[\left(\dfrac{ X_{j}}{b(x)}\,,\; \dfrac{e^{-R(\theta_1)}}{c(x)}\right) \in (0,\,\infty]^2\right]=\mu_j^*((0,\,\infty]^2)>0\,.
\eeam
Then we obtain	$H \in (\mathcal{D}\cap \mathcal{A})_A$, and if $\phi(p) < 0$, for some $p > \alpha  \beta/(\alpha+\beta)$, then this means that $\phi(p)< 0$ for some $p > J_{H_A}^+$.
\eexam
	
\pr~
From relations \eqref{eq.KLP.5.15} and  \eqref{eq.KLP.5.16}, we can apply \cite[Th. 4]{fougeres:mercadier:2012}, and thus we obtain that ${\bf X}\,e^{-R(\theta_1)}$ follows distribution 
\beao
H \in MRV\left(\dfrac {\alpha\,\beta}{\alpha + \beta}\,,\;\widehat{\mu} \right)\,,
\eeao
with the new measure $\widehat{\mu}$ 	to be such that $\widehat{\mu}(\bbb)= \widehat{\mu}\left[\left(y,\,{\bf x} \right)\;:\; y\,{\bf x} \in \bbb\right]$. Hence, from \eqref{eq.KLP.2.10} we find that $H \in (\mathcal{D}\cap \mathcal{A})_A$. Further, since 
\beao
J_{H_A}^+= \dfrac {\alpha\,\beta}{\alpha + \beta}\,,
\eeao 
we have the desired inequalities for the Laplace exponent.
~\halmos
	
\bre \label{rem.KLP.5.3}
The requirement of relation \eqref{eq.KLP.5.16}, means that at least one component of vector ${\bf X}$ is asymptotically dependent with the financial risk $e^{-R(\theta_1)}$. This is enough, to say that $X_A$ and $e^{-R(\theta_1)}$ are asymptotically dependent for any $A \in \mathscr{R}$. Indeed, since ${\bf X}$ has non-negative components and choosing as $j$, this of relation \eqref{eq.KLP.5.16}, and under the condition that $p_j>0$, we obtain that
\beao
&&\lim \dfrac{\PP[X_A > b(x)\,,\;e^{-R(\theta_1)}> c(x)]}{\PP[X_A > b(x)]} \geq \lim \dfrac{\PP\left[\dfrac{p_j\,X_j}{b(x)} > 1\,,\;\dfrac{e^{-R(\theta_1)}}{ c(x)} >1 \right]}{\PP\left[\dfrac{{\bf X}}{ b(x)} \in A \right]}\\[2mm]
&&=\dfrac{\mu_j^{*}((p_j,\,\infty]\times(1,\,\infty])}{\mu(A)}>0\,,
\eeao  
where at the third step we used relation \eqref{eq.KLP.5.16} and the fact that $\mu(A) \in (0,\,\infty)$, for any $A \in \mathscr{R}$, see in proof of \cite[Prop. 4.14]{samorodnitsky:sun:2016}. Quit similarly, when $p_j>0$, we find that it holds
\beao
&&\lim \dfrac{\PP[X_A > b(x)\,,\;e^{-R(\theta_1)}> c(x)]}{\PP\left[e^{-R(\theta_1)}> c(x)\right]}>0\,.
\eeao
Hence, when $p_j>0$, under the conditions of Example \ref{exam.KLP.5.2}, the $X_A$ and  $e^{-R(\theta_1)}$ are asymptotically dependent.
\ere
	
\bre \label{rem.KLP.5.4}
Under the conditions of Example \ref{exam.KLP.5.2}, we obtain that the insurance and financial risks, with regular variation parameters $\alpha,\,\beta$ respectively, have lighter tail from the final discounted claim, since $H_A \in \mathcal{R}_{-(\alpha\,\beta)/(\alpha + \beta)}$ and 
\beao
\dfrac{\alpha\,\beta}{\alpha + \beta} <\alpha \wedge \beta\,.
\eeao
We observe, that since 
\beao
H \in MRV\left(\dfrac{\alpha\,\beta}{\alpha + \beta}\,,\;\widehat{\mu} \right)\,,
\eeao
there exists a auxiliary distribution $G^* \in \mathcal{R}_{-(\alpha\,\beta)/(\alpha + \beta)}$, such that it holds
\beam \label{eq.KLP.5.18} 
\lim \dfrac 1{\bG^*(x)}\,\PP\left[{\bf X}\,e^{-R(\theta_1)}\in x\,\bbb \right]=\widehat{\mu}(\bbb)\,.
\eeam
\ere
	
\bco \label{cor.KLP.5.2}
Let $A \in \mathscr{R}$ be some fixed set, and let the discounted aggregate claims of relation \eqref{eq.KLP.1.3}. Under the Assumptions \ref{ass.KPT.2.1} and \ref{ass.KLP27.1.1} and with additional assumptions of Example \ref{exam.KLP.5.2} and if $\phi(p)<0$, 
for some 
\beao
p> \dfrac{\alpha\beta}{\alpha+\beta}\,,
\eeao 
then we obtain that
\beam \label{eq.KLP.5.19} \notag
&&\PP\left[{\bf D}(\infty)\in x\,A\right]\sim \widehat{\mu}(A) \,\bG^*(x)\left(1 + \int_{\Lambda} \int_0^{\infty} \E\left[e^{-\alpha\,\beta\,R(s)/(\alpha + \beta) } \right]\,\lambda(ds)\,\PP[\theta_1 \in du] \right)\,.\\
\eeam
\eco
	
\pr~
From Example \ref{exam.KLP.5.2} we obtain that Theorem \ref{th.KLP.4.1} holds. It remains to show that the right hand side in relation 	\eqref{eq.KLP.4.1} is asymptotic equivalent to the right hand side of relation \eqref{eq.KLP.5.19}. We mention 	that with the help of relation \eqref{eq.KLP.5.3}, for any $i \in \bbn$, we obtain that
\beam \label{eq.KLP.5.21} 
\E\left[e^{-p\,R(\tau_{i-1})} \right]=\E\left[e^{-(i-1)\,p\,R(\theta_1)} \right] < 1\,
\eeam 
for some $p >\alpha\,\beta/(\alpha + \beta)$. Hence, we obtain that
\beao
&&\int_0^{\infty} \PP\left[{\bf X}\, e^{-R(s)} \in x\,A \right]\,\lambda(ds)=\sum_{i=1}^{\infty}\PP\left[{\bf X}^{(i)}\, e^{-R(\tau_i)} \in x\,A\,\right] \\[2mm]
&&=\sum_{i=1}^{\infty}\PP\left[X_A^{(i)}\, e^{R(\tau_{i-1})-R(\tau_i)}\, e^{-R(\tau_{i-1})}> x \right]\\[2mm]
&&\sim \sum_{i=1}^{\infty} \E\left[ e^{-\alpha\,\beta\,R(\tau_{i-1})/(\alpha + \beta) } \right]\,\PP\left[X_A^{(i)}\, e^{R(\tau_{i-1})-R(\tau_i)}>x\,\right]\\[2mm] 
&&= \PP\left[X_A\, e^{-R(\theta_{1})}> x \right]\,\sum_{j=0}^{\infty} \E\left[ e^{-\alpha\,\beta\,R(\tau_{j})/(\alpha + \beta) } \,\right]\\[2mm]
&&\sim \widehat{\mu}(A)\,\bG^*(x)\,\left( 1 + \int_{\Lambda} \int_{0}^{\infty} \E\left[ e^{-\alpha\,\beta\,R(s)/(\alpha + \beta) } \right]\,\lambda(ds)\,\PP[\theta_1 \in du] \right)\,, 
\eeao
where at the third step we used Breiman theorem, in view of relation \eqref{eq.KLP.5.21}, at the fourth step we employ change of the variables in summation, while in last step we used the fact that 
\beao
H \in MRV\left(\dfrac{\alpha\,\beta}{\alpha + \beta},\,\widehat{\mu}\right)
\eeao
and relation \eqref{eq.KLP.5.18}. 
~\halmos

\subsection{Proof of Theorem \ref{th.KLP.4.1}}\label{sec.KLP.4.2}

Before the proof of the Theorem \ref{th.KLP.4.1} we need some preliminary lemmas. In what follows the $\{W_i\,,\;i \in \bbn \}$ denote a sequence of i.i.d. non-negative and non-degenerate to zero, random variables, and we write 
\beao
\Pi_i = \prod_{j=1}^i W_j\,.
\eeao 
We note that in the following lemma, we say that $\{X_i\,\Pi_i \,,\;i \in \bbn\}$ is tail asymptotically independent on $A$, symbolically $TAI_A$, for  $A \in \mathscr{R}$, when the $\{X_A^{(i)}\,\Pi_i \,,\;i \in \bbn \}$ is tail asymptotically independent, $TAI$, namely for any $i \neq j$, with $i,\,j \in \bbn$ we have that
\beam \label{eq.KLP.3.3}
\lim_{x_i \wedge x_j \to \infty} \PP[X_A^{(i)}\,\Pi_i > x_i\;|\;X_A^{(j)}\,\Pi_j > x_j]=0\,.
\eeam 
	
The TAI structure was introduced in \cite{geluk:tang:2009}, and was used to later papers, since  represents a rather general dependence structure, see for example in \cite{li:2013}, \cite{cheng:2014} and \cite{zou:peng:xu:2025}.
	
\ble \label{lem.KLP.3.1}
Let $A \in \mathscr{R}$ be a fixed set and let $\{{\bf X}^{(i)}\,,\;i \in \bbn\}$ be a sequence of i.i.d. non-negative random vectors with common distribution $F$. Let $\{{\bf Z}^{(i)}={\bf X}^{(i)}\,W_i \,,\;i \in \bbn\}$ be a sequence of i.i.d. random vectors with common distribution $V \in(\mathcal{D}\cap \mathcal{L})_A $. Let suppose that $\E[W^p] < \infty$, for some $p > J_{V_A}^+$. Then we have that
\begin{enumerate}
\item[(i)]
the $\{{\bf X}^{(i)}\,\Pi_i \,,\;i \in \bbn\}$ have distribution from class $(\mathcal{D}\cap \mathcal{L})_A $, and for any $i\neq j,\,i,j \in \bbn$ it holds 
\beao
\PP\left[{\bf X}^{(i)}\,\Pi_i \in x\,A\right] \asymp \PP\left[{\bf X}^{(j)}\,\Pi_j \in x\,A\right]\,.
\eeao
\item[(ii)]
The sequence $\{{\bf X}^{(i)}\,\Pi_i \,,\;i \in \bbn\}$ is $TAI_A$.
\end{enumerate} 
\ele
	
\pr~
\begin{enumerate}
\item[(i)]
We choose some ${\bf X}^{(i)}\,\Pi_i$, for $i \in \bbn$, and we show that it has distribution from class $(\mathcal{D}\cap \mathcal{L})_A $. For this, it is enough to show that the distribution of the product $X_A^{(i)}\,\Pi_i$ belongs to the class $\mathcal{D}\cap \mathcal{L}$, since we have 
\beao
\PP\left[{\bf X}^{(i)}\,\Pi_i \in x\,A\right] = \PP\left[X_A^{(i)}\,\Pi_i > x\right] \,,
\eeao
remember relation \eqref{eq.KLP.3.1b}. Therefore, since $X_A^{(i)}\,\Pi_i =X_A^{(i)}\,W_i\,\cdots\,W_1$, from assumption the product $X_A^{(i)}\,W_i$ has distribution $V_A \in \mathcal{D}\cap \mathcal{L}$, by $\E[W_{i-1}^p]< \infty$, for some $p > J_{V_A}^+$, and also by independence between $X_A^{(i)}\,W_i$ and $W_{i-1}$, from \cite[Th. 3.3(iv)]{cline:samorodnitsky:1994}, see also in \cite[Cor. 5.2]{leipus:siaulys:konstantinides:2023}, we obtain that the distribution of $X_A^{(i)}\,W_i\,W_{i-1}$ belongs to $\mathcal{D}\cap \mathcal{L}$ and further we have that
\beao
\PP\left[X_A^{(i)}\,W_i\,W_{i-1} > x\right] \asymp \PP\left[X_A^{(i)}\,W_i > x\right] \,,
\eeao 
hence 
\beao
\PP\left[{\bf X}^{(i)}\,\,W_i\,W_{i-1}  \in x\,A\right] \asymp \PP\left[X_A^{(i)}\,W_i \in x\,A \right] \,.
\eeao
We observe that from \cite[Lem. 3.9]{tang:tsitsiashvili:2003a} we find that the upper Matuszewska indexes of the distributions of $X_A^{(i)}\,W_i$ and $X_A^{(i)}\,W_i\,W_{i-1}$ are equal. So, continuing similarly for the $X_A^{(i)}\,W_i\,W_{i-1}\,W_{i-2} $ and so on, we conclude that the assertion (1) is true.
\item[(ii)]
It is enough to show that relation \eqref{eq.KLP.3.3} is true. Without loss of generality we can choose $i<j$, with $i,\,j \in \bbn$. Then we obtain that
\beao
&&\lim_{x_i \wedge x_j \to \infty}\PP\left[X_A^{(i)}\,\Pi_i > x_i\,,\;X_A^{(j)}\,\Pi_j > x_j\right] \\[2mm] \notag
&&= \lim_{x_i \wedge x_j \to \infty}\PP\left[X_A^{(i)}\,\Pi_i > x_i\,,\;X_A^{(j)}\,W_{j}\,\Pi_{j-1} > x_j\right] =o\left( \PP\left[X_A^{(j)}\,\Pi_j > x_j\right]\right)\,.
\eeao
where at the second step we apply \cite[Lem. 7]{tang:yuan:2014}. Indeed, as the $X_A^{(j)}\,W_{j}$ follows distribution $V_A$ belonging to class $\mathcal{D}\cap \mathcal{L} \subsetneq \mathcal{D}$, it is independent of $X_A^{(i)}\,\Pi_i$ and $\Pi_{j-1} $, since $i < j$, and $		\E[\Pi_{j-1}^p ]=\left(\E[W^p ]\right)^{j-1} < \infty$, for some $p> J_{V_A}^+$. For $i \geq j$ we apply the same argument in symmetric way.~\halmos
	\end{enumerate}
	
	The following lemma provides a type of multivariate linear single big-jump principle of $\widehat{\bf S}_{n}$, recall relation \eqref{eq.KPT.3.26}.
	
	\ble \label{lem.KLP.3.2}
	Under the conditions of Lemma \ref{lem.KLP.3.1}, for each $n \in \bbn$ we obtain that
	\beam \label{eq.KLP.3.8} 
	\PP\left[\widehat{\bf S}_{n} \in x\,A \right] \sim \sum_{i=1}^n\PP\left[{\bf X}^{(i)}\,\Pi_i \in x\,A\right]\,.
	\eeam
	\ele
	
	\pr~
	By Lemma \ref{lem.KLP.3.1}(i) we have that the products ${\bf X}^{(i)}\,\Pi_i$, for $i \in \bbn$, have distribution from class $(\mathcal{D}\cap \mathcal{L})_A$, hence the random variables $X_A^{(i)}\,\Pi_i \sim B_A^{(i)}$, where $B_A^{(i)} \in \mathcal{D}\cap \mathcal{L}$, for any $i =1,\,\ldots,\,n$. Hence, if $a_i$ is the insensitivity function for the distribution  $B_A^{(i)}$, then we have that
	\beam \label{eq.KLP.3.9} 
	a(x):= \bigwedge_{i=1}^n a_i(x)\,,
	\eeam
is the insensitivity function for all the distributions $B_A^{(1)},\,\ldots,\,B_A^{(n)}$. Hence, using at the first step \cite[Prop. 2.4]{konstantinides:passalidis:2024g} we obtain that
	\beam \label{eq.KLP.3.10} \notag
	&&\PP\left[\widehat{\bf S}_{n} \in x\,A \right] \leq \PP\left[\sum_{i=1}^n X_A^{(i)}\,\Pi_i > x\right]=\PP\Bigg[\sum_{i=1}^n  X_A^{(i)}\,\Pi_i > x\,,\;\bigvee_{j=1}^n X_A^{(j)}\,\Pi_j > x-a(x) \Bigg] \\[2mm] \notag
	&&+\PP\left[\sum_{i=1}^n  X_A^{(i)}\,\Pi_i > x\,,\;\bigvee_{j=1}^n X_A^{(j)}\,\Pi_j \leq x-a(x)\right]\leq \PP\left[\bigvee_{j=1}^n X_A^{(j)}\,\Pi_j > x-a(x) \right] \\[2mm] \notag
	&&+\PP\Bigg[\sum_{i=1}^n  X_A^{(i)}\,\Pi_i > x\,,\;\bigvee_{j=1}^n X_A^{(j)}\,\Pi_j \leq x-a(x)\,,\;\bigvee_{k=1}^n X_A^{(k)}\,\Pi_k >\dfrac xn\Bigg]\\[2mm] \notag
	&&\leq \sum_{i=1}^n \PP\left[X_A^{(i)}\,\Pi_i > x-a(x) \right] +\sum_{k=1}^n\PP\Bigg[\sum_{k\neq i=1}^n  X_A^{(i)}\,\Pi_i > a(x)\,,\; X_A^{(k)}\,\Pi_k >\dfrac xn\Bigg]\\[2mm] \notag
	&&\leq \sum_{i=1}^n \PP\left[X_A^{(i)}\,\Pi_i > x-a(x) \right] +\sum_{k=1}^n\sum_{k\neq i=1}^n\PP\Big[  X_A^{(i)}\,\Pi_i > \dfrac{a(x)}{n-1}\,,\; X_A^{(k)}\,\Pi_k >\dfrac xn\Big]\\[2mm] 
	&&\sim \sum_{i=1}^n \PP\left[X_A^{(i)}\,\Pi_i > x \right] =\sum_{i=1}^n\PP\left[  {\bf X}^{(i)}\,\Pi_i \in x\,A \right]\,,
	\eeam
where at the penultimate step we used the property of the insensitivity function in relation \eqref{eq.KLP.3.9} in first term, while any term of the second sum is negligible with respect to each term of the first one, since by Lemma \ref{lem.KLP.3.1}(ii) the sequence $\{X_A^{(i)}\,\Pi_i\,,\;i \in \bbn \}$ is $TAI$, and $B_A^{(i)} \in \mathcal{D}\cap \mathcal{L} $.
	
Hence, from relation \eqref{eq.KLP.3.10} we obtain the upper bound in \eqref{eq.KLP.3.8}. Further, we proceed to the lower bound. By the fact that the summands are non-negative, the set $A$ is increasing, and using the Bonferroni inequality we have that
\beam \label{eq.KLP.3.11} 
&&\PP\left[\widehat{\bf S}_{n} \in x\,A \right] \geq \PP\left[\bigcup_{i=1}^n \left\{ {\bf X}^{(i)}\,\Pi_i \in x\,A \right\} \right]\geq\\[2mm] \notag
&& \sum_{i=1}^n \PP\left[ {\bf X}^{(i)}\,\Pi_i \in x\,A \right]- \sum_{1\leq i < j \leq n} \PP\left[{\bf X}^{(i)}\,\Pi_i \in x\,A\,,\;{\bf X}^{(j)} \Pi_j \in x A \right]\sim \sum_{i=1}^n\PP\left[  {\bf X}^{(i)} \Pi_i \in x\,A \right] ,
\eeam
where at the last step we applied Lemma \ref{lem.KLP.3.1}(ii). Therefore, from relation \eqref{eq.KLP.3.11} we reach to the desired lower bound in \eqref{eq.KLP.3.8}.
~\halmos

In next Lemma we reformulate a result from \cite[Lem. 3.3]{chen:yuan:2017}, see also in \cite[Lem. 3.2]{hao:tang:2012}.

\ble \label{lem.KLP.4.1}
Let $Z$ be a real valued random variable with distribution $G$ and Matuszewska indexes $0 < J_{G}^- \leq J_G^+ < \infty$, that means $G \in \mathcal{D}\cap \mathcal{P_D}$. Then for any $0 < q_1 < J_{G}^- \leq J_G^+ < q_2 < \infty$, there exists some constant $C>0$ and some real number $x_0>0$, that depends on $q_1,\,q_2 $, such that for all $x> x_0$, and any non-negative random variable $\xi$, independent of $Z$, we have that
\beao
\dfrac{\PP\left[\xi\,Z >x\right]}{\PP\left[Z >x\right]} \leq C\,\E\left[\xi^{q_1}\vee \xi^{q_2} \right]\,.
\eeao
\ele

The following result plays a crucial role for the argumentation of Theorem \ref{th.KLP.4.1}.

\ble \label{lem.KLP.4.2}
Let $A \in \mathscr{R}$ be a fixed set and $\{{\bf X}^{(i)}\,,\; i \in \bbn\}$ be a sequence of i.i.d., non-negative random vectors with common distribution $F$. Let $\{{\bf Z}^{(i)}={\bf X}^{(i)}\,W_i\,,\;i \in \bbn\}$ be a sequence of i.i.d. random vectors with common distribution $V \in (\mathcal{D}\cap \mathcal{A})_A$, and Matuszewska indexes $0 < J_{V_A}^- \leq J_{V_A}^+ < \infty$. Let assume that there exists $q_1,\,q_2$, with $0 < q_1 < J_{V_A}^- \leq J_{V_A}^+ < q_2 < \infty$, such that it holds $\E\left[W^{q_1} \vee W^{q_2}\right]< 1$. Then \eqref{eq.KLP.3.8} holds uniformly for $n \in \bbn$, in the sense
\beam \label{eq.KLP.4.3} 
\lim_{x\to\infty} \sup_{n \in \bbn} \left|\dfrac{\PP\left[\widehat{\bf S}_{n} \in x\,A\right]}{\sum_{i=1}^n\PP\left[{\bf X}^{(i)}\,\Pi_i \in x\,A\right]} -1 \right| =0\,.
\eeam
\ele

\pr~
We will follow the steps of the proof of Theorem 3.3 of \cite{chen:yuan:2017}, under some modifications. Let $M \in \bbn$ some arbitrarily large integer number. Then by Lemma \ref{lem.KLP.3.2} we obtain that relation \eqref{eq.KLP.4.3} holds uniformly for $n=1,\,\ldots,\,M$. Thus, we proceed with showing that relation \eqref{eq.KLP.4.3} holds uniformly for $n > M$. Initially, from \cite[Eq. (3.4)]{chen:yuan:2017} we obtain that 
\beam \label{eq.KLP.4.4} 
\lim_{M \to \infty} \limsup \dfrac{\sum_{i=M+1}^{\infty}\PP\left[{\bf X}^{(i)}\,\Pi_i \in x\,A\right]}{\PP\left[{\bf X} W \in x A\right]}  =\lim_{M \to \infty} \limsup \dfrac{\sum_{i=M+1}^{\infty}\PP\left[X_A^{(i)} \Pi_i > x\right]}{\PP\left[X_A W > x\right]} =0 ,
\eeam
where the last step was shown in \cite{chen:yuan:2017} for the case $V_A \in \mathcal{C}\cap \mathcal{P_D}$, however the same argument still holds for the case $V_A \in \mathcal{D}\cap \mathcal{P_D}$, see Lemma \ref{lem.KLP.4.1}, and consequently for the case $V_A \in \mathcal{D}\cap \mathcal{A}$. Further, similarly by \cite[Eq. (3.5)]{chen:yuan:2017} in third step we obtain that
\beam \label{eq.KLP.4.5}
&&\limsup \dfrac{\PP\left[\sum_{i=1}^{\infty}{\bf X}^{(i)}\,\Pi_i \in x\,A\right]}{\PP\left[{\bf X}\,W \in x\,A\right]}  = \limsup \dfrac{\PP\left[ \sup_{{\bf p} \in I_A} {\bf p}^T\,\sum_{i=1}^{\infty}{\bf X}^{(i)}\,\Pi_i > x\right]}{\PP\left[X_A\,W > x\right]} \\[2mm]\notag
&& \leq \limsup \dfrac{\PP\left[\sum_{i=1}^{\infty}X_A^{(i)}\,\Pi_i > x\right]}{\PP\left[X_A\,W > x\right]} < \infty\,.
\eeam
Let us introduce a new random variable $\eta$, independent from all the other sources of randomness, that satisfies the relations
\beam \label{eq.KLP.4.6} 
\PP\left[\eta > x\right]\sim c^*\,\PP\left[X_A\,W > x\right]\,,
\eeam
for some constant $c^*>0$, and 
\beam \label{eq.KLP.4.7} 
\PP\left[\sum_{i=1}^{\infty}X_A^{(i)}\,\Pi_i > y\right] \leq \PP\left[\eta > y\right]\,,
\eeam
for any $y \in \bbr$ (note that the existence \eqref{eq.KLP.4.7} is ensured by \eqref{eq.KLP.4.5}). For the upper bound in \eqref{eq.KLP.4.3}, for $n>M$ we have that
\beam \label{eq.KLP.4.8} \notag
&&\PP\left[\widehat{\bf S}_{n} \in x\,A\right] \leq \PP\left[ \sum_{i=1}^{n}X_A^{(i)}\,\Pi_i > x\right] \\[2mm]\notag
&& =\PP\left[\sum_{i=1}^{M}X_A^{(i)}\,\Pi_i +\left(\sum_{i=M+1}^{n}X_A^{(i)}\,\prod_{j=M+1}^i W_j \right)\Pi_M> x\right]\\[2mm]
&& \leq \PP\left[\sum_{i=1}^{M}X_A^{(i)}\,\Pi_i +\eta\,\Pi_M> x\right]\sim
\sum_{i=1}^{M}\PP\left[X_A^{(i)}\,\Pi_i > x\right] + \PP\left[\eta\,\Pi_M> x\right]\,,
\eeam
where at the first step we applied \cite[Prop. 2.4]{konstantinides:passalidis:2024g}, at the third step we used relation \eqref{eq.KLP.4.7}. At the fourth step we applied \cite[Th. 3.1]{geluk:tang:2009}, since the products $X_A^{(1)}\,\Pi_1\,,\;\ldots\,,\;X_A^{(M)}\,\Pi_M,\eta\,\Pi_M$ have distributions in class $\mathcal{D}\cap \mathcal{L}$ and are $TAI$. Indeed, from Lemma \ref{lem.KLP.3.1}(i), the products $X_A^{(1)}\,\Pi_1\,,\;\ldots\,,\;X_A^{(M)}\,\Pi_M$ have distributions that belong to $\mathcal{D}\cap \mathcal{L}$, while by Lemma \ref{lem.KLP.3.1}(ii) these products are 
$TAI$. Further, from relation \eqref{eq.KLP.4.6}, in view of closure property of class $\mathcal{D}\cap \mathcal{L}$ with respect to strong tail equivalence, see in \cite[p. 80]{leipus:siaulys:konstantinides:2023}, we obtain 
that $\eta$ has distribution belonging to $\mathcal{D}\cap \mathcal{L}$ (more precisely belongs to $\mathcal{D}\cap \mathcal{A}$), and since there exists some $p$ such that $J_{V_A}^+ < p < \infty$, and the inequality $\E\left[\Pi_M^{p}\right]= \left(\E\left[W^{p}\right]\right)^M< 1 < \infty$ holds, we find that
\beao
\PP[\Pi_M > v\,x] =o\left(\PP\left[X_A \,W> x\right]\right) =o\left(\PP\left[\eta> x\right]\right)\,, 
\eeao
for any $v>0$. Hence, by \cite[Th. 2.2(iii), 3.3(ii)]{cline:samorodnitsky:1994} we have that the product $\eta\,\Pi_M$ has distribution from class $\mathcal{D}\cap \mathcal{L}$. Further we choose some $X_A^{(i)}\,\Pi_i$, for some $i=1,\,\ldots,\,M$ to obtain that
\beao
&&\PP[\eta\,\Pi_M > x_j\,,\; X_A^{(i)}\,\Pi_i>x_i]=o\left(\PP\left[X_A^{(i)} \,\Pi_i> x_i \right]\right)\,,\\[2mm]
&&\PP[\eta\,\Pi_M > x_j\,,\; X_A^{(i)}\,\Pi_i>x_i]=o\left(\PP\left[\eta\,\Pi_M> x_j \right]\right)\,,
\eeao
as $x_i \wedge x_j \to \infty$, where we applied \cite[Lem. 7]{tang:yuan:2014},  in case of $i=M$, the first relation still holds, but this time showing that the $\eta$ and $X_A\,W$ are $TAI$ from \cite[Lem. 7]{tang:yuan:2014}, and further the $\eta\,\Pi_M$ and $X_A\,\Pi_M$, continue to be $TAI$ by \cite[Th. 2.2]{li:2013} (recall that from \eqref{eq.KLP.4.6}, the $\eta$ has upper Matuszewska index $J_{V_A}^+$). Hence, we can apply \cite[Th. 3.1]{geluk:tang:2009} at the fourth step of \eqref{eq.KLP.4.8}.

For any $\vep >0$, from relation \eqref{eq.KLP.4.6} and Lemma \ref{lem.KLP.4.1}, via moment conditions, we can find large enough $M$ and some constant $C>0$ such that we have
\beam \label{eq.KLP.4.9} \notag
&&\dfrac{\PP\left[\eta\,\Pi_M > x\right]}{\PP\left[X_A \,W> x\right]} \sim c^*\,\dfrac{\PP\left[\eta\,\Pi_M > x\right]}{\PP\left[\eta > x\right]} \leq c^*\,C\,\E\left[\Pi_M^{q_1} \vee \Pi_M^{q_2}\right]\\[2mm]
&& \leq c^*\,C\,\left(\E\left[\left(W^{q_1}\right)^M\right] + E\left[ \left(W^{q_2}\right)^M\right] \right) \leq \vep\,.
\eeam  
From \eqref{eq.KLP.4.8} and \eqref{eq.KLP.4.9} we find that for any  $n>M$ it holds
\beam \label{eq.KLP.4.10}
\PP\left[\widehat{\bf S}_{n} \in x\,A\right] \lesssim (1+\vep)\,\sum_{i=1}^{M}\PP\left[X_A^{(i)} \,\Pi_i > x\right] \leq (1+\vep)\,\sum_{i=1}^{n}\PP\left[{\bf X} \,\Pi_i  \in x\,A \right]\,.
\eeam
From the other hand, by Lemma \ref{lem.KLP.3.2}, for any $n> M$ we obtain that
\beam \label{eq.KLP.4.11}
&&\PP\left[\widehat{\bf S}_{n} \in x\,A\right] \geq \PP\left[\sum_{i=1}^{M}{\bf X}^{(i)}\,\Pi_i \in x\,A\right] \sim \sum_{i=1}^{M}\PP\left[{\bf X}^{(i)}\,\Pi_i \in x\,A\right]\\[2mm] \notag
&&= \left(\sum_{i=1}^{n} - \sum_{i=M+1}^{n}\right) \PP\left[{\bf X}^{(i)}\,\Pi_i \in x\,A\right] \gtrsim (1-\vep)\,\sum_{i=1}^{n} \PP\left[{\bf X}^{(i)}\,\Pi_i \in x\,A\right]\,,
\eeam
where at the last step we used relation \eqref{eq.KLP.4.4}. Therefore, by relations \eqref{eq.KLP.4.10}, \eqref{eq.KLP.4.11} and the arbitrariness in choice of $\vep>0$, we obtain \eqref{eq.KLP.4.3} uniformly for $n> M$.
~\halmos

{\bf Proof of Theorem \ref{th.KLP.4.1}}
Since $\phi(p_2)< 0$, for some $p_2>J_{H_A}^+$, that means $\phi(z)< 0$ for any $z \in (0,\,p_2]$, which follows from the facts that $\phi(0)= 0$ and $\phi(z)$ is convex function with respect to $z$. Hence, from relation \eqref{eq.KLP.1.1} we obtain that
\beao
\E\left[e^{-p_2\,R(\theta_1)} \right]=\E\left[e^{\theta_1\,\phi(p_2)} \right]<1\,,
\eeao
and further for $0 < p_1 < J_{H_A}^- \leq J_{H_A}^+ < p_2 < \infty$, we have that
\beam \label{eq.KLP.4.12}
\E\left[e^{-p_1\,R(\theta_1)} \vee e^{-p_2\,R(\theta_1)}  \right]< 1\,.
\eeam
From relation \eqref{eq.KLP.4.12}, in combination with the relation 
\beao
W_i = e^{R(\tau_{i-1}) -R(\tau_i)} \stackrel{d}{=} e^{-R(\theta_1)}\,,
\eeao 
for any $i \in \bbn$, applying Lemma \ref{lem.KLP.4.2} for $n=\infty$, we obtain that
\beao
&&\PP\left[{\bf D}(\infty) \in x\,A\right] = \PP\left[\sum_{i=1}^{\infty} {\bf X}^{(i)}\,e^{-R(\tau_i)} \in x\,A \right]\\[2mm]
&&\sim \sum_{i=1}^{\infty} \PP\left[ {\bf X}^{(i)}\,e^{-R(\tau_i)} \in x\,A \right] =\int_0^{\infty} \PP\left[ {\bf X}\,e^{-R(s)} \in x\,A \right]\,\lambda(ds)\,,
\eeao
which completes the proof.
~\halmos

\section{Applications on ruin probabilities} \label{sec.KPT.4}

In risk theory, when the risk models have claim distributions with heavy tails, 
then the asymptotic behavior of the ruin probability coincides with the 'tail' asymptotic behavior of the discounted aggregate claims. There are many papers on this direction on one-dimensional or multidimensional set up, see for example in \cite{tang:wang:yuen:2010}, \cite{li:2012}, \cite{chen:konstantinides:passalidis:2025}, \cite{konstantinides:passalidis:2024j}, among others.

In this section we study the insurer's ruin probability, having in mind the asymptotic 
behavior of the discounted aggregate claims, as was analyzed in Sections 3 and 4. We consider that the insurer's discounted  surplus process at the moment $t\geq 0$ is given by the relation
\beam \label{eq.KPT.4.1} 
{\bf U}(t)=\left( 
\begin{array}{c}
U_{1}(t) \\ 
\vdots \\ 
U_{d}(t)
\end{array} 
\right) = x\,{\bf l} + \int_0^t e^{-R(s)}\,{\bf C}(ds) - {\bf D}(t)\,,
\eeam
where $x>0$ is the insurer's initial capital, that is allocated through the vector ${\bf l}=(l_1,\,\ldots,\,l_d)$, to $d$-lines of business, namely $l_1,\,\ldots,\,l_d >0$ with $l_1+\cdots+l_d =1$. The middle term in the right member of relation \eqref{eq.KPT.4.1} is understood as
\beao 
\int_0^t e^{-R(s)}\,{\bf C}(ds)=\left( 
\begin{array}{c}
\int_0^t e^{-R(s)}\,C_1(ds) \\ 
\vdots \\ 
\int_0^t e^{-R(s)}\,C_d(ds)
\end{array} 
\right) \,,
\eeao
where $\{{\bf C}(t) = (C_1(t),\,\ldots,\,C_d(t))\,,\;t\geq 0\}$ are the premiums of the $d$-lines of business, that represent non-negative c\'{a}dl\'{a}g processes with ${\bf C}({\bf 0}) ={\bf 0}$.  

In the multivariate risk model, the ruin probability is defined in several forms, we refer to \cite{cheng:yu:2019} and \cite{lu:li:yuan:shen:2024}, for more discussions. We give a general definition of ruin probability in multivariate risk model, as the 
entrance probability of the surplus process into some 'ruin' set $L$. Here, we use the 
ruin set $L$, given by the following assumption. 

\begin{assumption} \label{ass.KPT.4.1}
Let $L$ be a set, that is open, decreasing, that means the $-L$ is increasing, and it satisfies the conditions ${\bf 0} \in \partial L$, $L^c$ is convex, and $x\,L=L$ for any $x>0$. 
\end{assumption}

\bre \label{rem.KPT.4.1}
The ruin set of Assumption \ref{ass.KPT.4.1} was introduced in \cite{samorodnitsky:sun:2016}, while we easily find that the set 
\beam \label{eq.KPT.4.4} 
A:=({\bf l} -L) \in \mathscr{R}\,,
\eeam
when set $L$ satisfies Assumption \ref{ass.KPT.4.1}. From relation \eqref{eq.KPT.4.4} 
we can obtain some ruin sets, that are interesting for actuarial practice. First of all, the classical one-dimensional ruin set $L=(-\infty,\,0)$, which is related with the set $A = (1,\,\infty)$, that indicates the probability that the insurer's surplus becomes negative. 

Further, the set $A_2 \in \mathscr{R}$, from Remark \ref{rem.KPT.3.1}, that is related to the ruin set 
\beao
L_2 = \left\{{\bf y}\;:\;\sum_{j=1}^d y_j < 0 \right\}\,,
\eeao  
where the surplus entrance probability into $L_2$, is the probability that the sum of 
surpluses of insurer's $d$-lines of business, to become some moment negative. Additionally, another important ruin set, that is related with set $A_1 \in \mathscr{R}$ of Remark \ref{rem.KPT.3.1}, is the following $L_1 =\{ {\bf y}\;:\;y_i < 0\,,\;\exists \; i=1,\,\ldots,\,d\}$, where the entrance probability of the surplus in $L_1$, indicates the probability one of the insurer's $d$-lines of business to have some moment negative surplus. For more examples and discussions about the ruin sets of Assumption \ref{ass.KPT.4.1} we refer to \cite[Sec. 5]{samorodnitsky:sun:2016} and \cite[Sec. 2]{shen:yuan:lu:2022}.
\ere

So, we can define the ruin probability based on set $L$, which in what follows will satisfy Assumption \ref{ass.KPT.4.1}. The infinite time ruin probability, is defined, via  Assumption \ref{ass.KPT.4.1} and 
relation \eqref{eq.KPT.4.4}, as follows
\beam \label{eq.KPT.4.5} 
&&\psi_{{\bf l},\,L}(x\,;\;\infty):=\PP\left( {\bf U}(t) \in L\,,\;\exists \; t \geq 0 \right)\\[2mm] \notag
&& =\PP\left( {\bf D}(t) -\int_0^t e^{-R(s)}\,{\bf C}(ds) \in x\,({\bf l} - L)\,,\;\exists \;  t \geq 0 \right)\\[2mm] \notag
&& =\PP\left( {\bf D}(t) -\int_0^t e^{-R(s)}\,{\bf C}(ds) \in x\,A\,,\;\exists \; t \geq 0 \right)\,.
\eeam 

We need the technical assumption, that the  discounted aggregate premiums of each line of business $i=1,\,\ldots,\,d$ is such that it holds
\beam \label{eq.KPT.A}
0 \leq \widehat{C}_{i}(\infty):= \int_0^{\infty} e^{-R(s)}\,C_{i}(ds)<\infty\,,
\eeam
almost surely. We also define
\beam \label{eq.KPT.5.10} 
\widehat{C}(\infty):= \bigvee_{i=1}^d \widehat{C}_{i}(\infty)\,.
\eeam

The following result provides the infinite time ruin probability, under the conditions of Theorem \ref{th.KPT.3.1} and Theorem \ref{th.KLP.4.1}.
\bco \label{cor.KPT.5.2} 
Let $A = ({\bf l} - L) \in \mathscr{R}$ be some fixed set, and let consider the risk model in equation \eqref{eq.KPT.4.1} with premiums that satisfy relation \eqref{eq.KPT.A}. If the premiums are independent from all other sources of randomness, then:
\begin{enumerate}
\item[(i)]
If the conditions of Theorem \ref{th.KPT.3.1} hold, then we obtain
\beam \label{eq.KPT.5.11} 
\psi_{{\bf l},\,L}(x\,;\;\infty) \sim \int_0^{\infty} \PP\left( {\bf X}\,e^{-R(s)} \in x\,A \right)\,\lambda(ds)\,.
\eeam
\item[(ii)]
If the conditions of Theorem \ref{th.KLP.4.1} hold, then relation \eqref{eq.KPT.5.11} holds.
\end{enumerate}
\eco 

\pr~
We show the assertion (ii), since the (i) follows with similar steps. For the bound from above of the infinite horizon ruin probability, we take into account relation \eqref{eq.KPT.4.5}, the fact that the set $A$ is increasing and the $ {\bf D}(t)$ is non decreasing function of $t$, to find that
\beao
\psi_{{\bf l},\,L}(x\,;\;\infty)\leq \PP[ {\bf D}(\infty) \in x\,A] \sim \int_0^{\infty} \PP\left( {\bf X}\,e^{-R(s)} \in x\,A \right)\,\lambda(ds)\,,
\eeao
where at the second step we used Theorem \ref{th.KLP.4.1}. Further, we study the lower bound of relation \eqref{eq.KPT.5.11}, recalling that the premiums are independent of the other sources of randomness. We mention that by Theorem \ref{th.KLP.4.1}, or even by Lemma \ref{lem.KLP.4.2}, we obtain the asymptotic relation
\beam \label{eq.KPT.5.12} 
\PP\left( {\bf D}(\infty) \in x\,A \right) \sim \sum_{i=1}^{\infty} \PP\left( {\bf X}^{(i)}\,e^{-R(\tau_i)} \in x\,A \right)\,,
\eeam
while by Lemma \ref{lem.KLP.3.1}(i) we find that each product ${\bf X}^{(i)}\,e^{-R(\tau_i)} $ belongs to class $(\mathcal{D} \cap \mathcal{L})_A$. Hence, in view of relations \eqref{eq.KPT.4.5} and \eqref{eq.KPT.5.10} we obtain that
\beam \label{eq.KPT.5.13} \notag
&&\psi_{{\bf l},\,L}(x\,;\;\infty)\geq \PP\left[ {\bf D}(\infty)- \widehat{C}(\infty) \in x\,A\right]= \int_0^{\infty}\PP\left( {\bf D}(\infty) \in x\,A+y \right) \,\PP\left( \widehat{C}(\infty) \in dy \right)\\[2mm] \notag
&& \geq \int_0^{\infty}\PP\left( {\bf D}(\infty) \in (x+k)\,A \right) \,\PP\left( \widehat{C}(\infty) \in dy \right)\sim \sum_{i=1}^{\infty} \PP\left( {\bf X}^{(i)}\,e^{-R(\tau_i)} \in (x+k)\,A \right)\\[2mm]
&& \sim \sum_{i=1}^{\infty} \PP\left( {\bf X}^{(i)}\,e^{-R(\tau_i)} \in x\,A \right) = \int_0^{\infty}\PP\left( {\bf X}\,e^{-R(s)} \in x\,A \right) \,\lambda(ds) \,,
\eeam
where at the third step we used \cite[Lem. 4.3(d)]{samorodnitsky:sun:2016} for some  $0<k<x$, at the fourth step we employed relation \eqref{eq.KPT.5.12}, while at the fifth step we recalled the properties of distribution class $\mathcal{L}$ for the product $ X_A^{(i)}\,e^{-R(\tau_i)} $. Thus, by relation \eqref{eq.KPT.5.13} we get the lower bound. We note that for the proof of assertion (i), we need the property of class $\mathcal{L}_A$, for the products ${\bf X}^{(i)}\,e^{-R(\tau_i)} $. Firstly we have that the products ${\bf X}^{(i)}\,e^{-R(\tau_i)}\stackrel{d}{=} {\bf X}^{(i)}\,e^{-R(\tau_{i-1})}\,e^{-R(\theta_1)}$. From \eqref{eq.KPT.2.11} we find that 
${\bf X}^{(i)}\,e^{-R(\theta_1)} \in \mathcal{A}_A$, while $e^{-R(\tau_{i-1})}\in [0,\,1]$, and is independent of ${\bf X}^{(i)}\,e^{-R(\theta_1)}$. Therefore, by \cite[Theorem 4.1]{konstantinides:passalidis:2024h}, we find that 
${\bf X}^{(i)}\,e^{-R(\tau_i)} \in \mathcal{A}_A \subsetneq \mathcal{L}_A$.
~\halmos 

\bre \label{rem.KPT.4.3}
We observe that under the conditions of Corollaries \ref{cor.KLP.5.1} and \ref{cor.KLP.5.2}, the ruin probability in relation \eqref{eq.KPT.5.11} takes the following simple expressions
\beao
&&\psi_{{\bf l},\,L}(x\,;\;\infty)\sim \mu(A)\,\bG(x)\,\\[2mm]
&&\times \E\left[e^{-\alpha\,R(\theta_1)}\,h\left(e^{-R(\theta_1)} \right) \right] \left(1 + \int_{\Lambda} \int_0^{\infty} \E\left[e^{-\alpha\,R(s)} \right]\,\lambda(ds)\,\PP(\theta_1 \in du) \right)\,,
\eeao
and
\beao
\psi_{{\bf l},\,L}(x\,;\;\infty)\sim \widehat{\mu}(A)\,\bG^*(x)\,\left(1 + \int_{\Lambda} \int_0^{\infty} \E\left[e^{-\alpha\,\beta\,R(s)/(\alpha + \beta)} \right]\,\lambda(ds)\,\PP(\theta_1 \in du) \right)\,,
\eeao
respectively.
\ere

\section*{Acknowledgments}
We feel the pleasant duty to express our sincere gratitude to Zhangting Chen, for several discussions, which substantially improved the paper.

\end{document}